\documentclass[journal]{IEEEtran}
\usepackage{caption}
\usepackage{subcaption}
\usepackage{graphicx}
\usepackage{algorithm}
\usepackage{algorithmicx}
\usepackage{algpseudocode}
\usepackage{makecell}
\usepackage{amsmath} 
\usepackage{amssymb}  
\usepackage{diagbox}
\usepackage{color}
\usepackage{multirow}

\graphicspath{{figures/}}

\newcommand{\etal}{\textit{et al}.~}
\newcommand{\ie}{\textit{i}.\textit{e}.}
\newcommand{\eg}{\textit{e}.\textit{g}.}

\newtheorem{definition}{Definition}
\newtheorem{remark}{Remark}
\newtheorem{lemma}{Lemma}
\newtheorem{theorem}{Theorem}

\newtheorem{assumption}{Assumption}

\ifCLASSINFOpdf
\else
\fi
\begin{document}
%
\title{Cooperation Method of Connected and Automated Vehicles at Unsignalized Intersections:\\ Lane Changing and Arrival Scheduling}
%
%
%

\author{
	Chaoyi Chen$^{1}$,
	Mengchi Cai$^{1}$,
	Jiawei Wang$^{1}$,
	Kai Li$^{2}$,
	Qing Xu$^{1}$, 
	Jianqiang Wang$^{1}$, 
	Keqiang Li*$^{1}$.
	\thanks{*This work was supported by the National Key Research and Development Program of China under Grant 2019YFB1600804, the National Natural Science Foundation of China under Grant 52072212, China Intelligent and Connected Vehicles (Beijing) Research Institute Co., Ltd. and Dongfeng Automobile Co., Ltd.}
	\thanks{$^{1}$Chaoyi Chen, Mengchi Cai, Qing Xu, Jianqiang Wang, and Keqiang Li are with the School of Vehicle and Mobility, Tsinghua University}
	\thanks{$^{2}$Kai Li is with Dongfeng Automobile Co., Ltd}
	\thanks{Corresponding author: Keqiang Li, Email address: likq@tsinghua.edu.cn}
}
\markboth{IEEE Transactions on Vehicular Technology,~Vol.~XX, No.~XX, XXX~2021}
{}
\maketitle

\begin{abstract}
The cooperation of connected and automated vehicles (CAVs) has shown great potential in improving traffic efficiency during intersection management. Existing research mainly focuses on intersections where lane changing is prohibited, which is impractical for real-life implementation. This paper proposes a two-stage cooperation framework, which decouples the longitudinal and lateral control of CAVs, allowing them to change to their preferred lanes. Based on formation control, an iterative framework is initially proposed to solve the target assignment and path planning problem for multiple CAVs on multiple lanes. A graph-based minimum clique cover method is then applied to obtain the optimal scheduling plan for the CAVs. Extensive numerical simulations for different numbers of vehicles and traffic volumes validate the effectiveness of the proposed algorithm.
\end{abstract}


%
\IEEEpeerreviewmaketitle

\section{Introduction}
Intersections are the most complicated scenarios in urban traffic, where traffic jams and vehicle collisions frequently occur~\cite{azimi2014stip}. With the rapid development of vehicle-to-everything (V2X) technology, the central coordinator deployed at an intersection can guide connected and automated vehicles (CAVs) through it. This guarantees the high-efficiency and conflict-free cooperation of CAVs~\cite{chen2020mixed}. After receiving the scheduled arrival time from the coordinator, each CAV optimizes its speed trajectory in pursuit of high traffic efficiency and low fuel consumption~\cite{xu2017v2i, guler2014using}. Several methods have been proposed to solve the longitudinal control problems of CAVs, including model predictive control~\cite{asadi2011predictive,yang2017eco}, fuzzy logic~\cite{milanes2009controller,onieva2015multi}, and optimal control~\cite{li2015eco,jiang2017eco}.


Apart from vehicle control, the CAV scheduling problem is also frequently discussed in research on intersections. CAVs approach from different directions with different destinations; thus, their trajectories inevitably intersect in the middle of the intersection. Therefore, staggering the arrival times of CAVs is an essential functionality of intersection management. Previous studies have also found that scheduling CAVs is the key factor influencing the traffic efficiency at intersections~\cite{li2006cooperative, meng2017analysis}. To schedule the CAVs approaching intersections, the most straightforward method is a first-in-first-out (FIFO) strategy, wherein the CAVs that enter first are scheduled to leave the intersection first~\cite{malikopoulos2018decentralized, dresner2008multiagent}. Similar concepts are applied in reservation-based~\cite{dresner2004multiagent}, batch-based~\cite{tachet2016revisiting}, and platoon-based methods~\cite{jin2013platoon}. 

Evidently, these ad hoc scheduling methods have low computational burden; however, they are less likely to obtain a high-efficiency scheduling plan. Therefore, other studies have also proposed optimization-based methods to solve the scheduling problem. Several researchers formulated this problem as a mixed integer program (MIP) problem~\cite{lu2018mixed, wang2020mixed, ge2021real}. Other methods, such as Monte Carlo tree search~\cite{xu2019cooperative}, dynamic programming~\cite{yan2009autonomous}, and minimum clique cover (MCC)~\cite{chen2021conflict}, were also proposed to schedule CAVs. It is worth mentioning that because scheduling CAVs is a discrete problem rather than a continuous one, a graph-based method is another promising method to solve this problem. Apart from the depth-first spanning tree algorithm proposed in~\cite{xu2018distributed}, a Petri net~\cite{lin2019graph} and a conflict duration graph~\cite{deng2020conflict} have also been used in modeling the scheduling problem. Although these studies widely investigated the CAV scheduling problem, the intersection scenarios remained limited to where lane changing was prohibited. 

In the aforementioned studies on intersections, the CAVs were assumed to run in their target lanes,~\ie, only their longitudinal control was considered. In practice, however, CAVs approach from random lanes and have different target lanes; therefore, it is necessary to extend this research to scenarios that permit lane changing. Earlier studies on this topic focused on obtaining a smooth CAV speed trajectory~\cite{li2005cooperative, li2007cooperative}; however, traffic efficiency was not investigated completely. With regard to CAV scheduling, which we are concerned with, a few studies formulated this problem as an MIP~\cite{lu2019trajectory} or a linear programming problem~\cite{hu2019trajectory}. This was done by assuming that lane changing maneuvers are accomplished in a given time interval. \cite{xu2020bi} proposed a practical bi-level framework, where the high-efficiency arrival plans and collision-free path planning are solved on the upper and lower levels separately. Several other prospective studies focused on changing lane directions dynamically rather than allocating CAVs to constant directional lanes~\cite{he2018erasing, cai2021multi, mitrovic2019combined}. However, flexible lane directions are unsuitable for a mixed traffic environment, where human-driven vehicles (HDVs) and CAVs coexist.

In summary, most of the existing studies were conducted in scenarios where lane changing was prohibited. Therefore, only CAV scheduling or vehicle longitudinal control were studied. In this study, we focus on a scenario where lane changing is permitted. CAVs are required to change to their target lanes on approaching the intersection. Most of the aforementioned studies circumvent this problem by assuming that the CAVs are running in their target lanes~\cite{li2006cooperative, malikopoulos2018decentralized, dresner2004multiagent, tachet2016revisiting, xu2018distributed}. A few studies considered lane changing behavior in their scheduling algorithm~\cite{lu2019trajectory,hu2019trajectory,xu2020bi}; however, the optimality of the scheduling was not extensively discussed. Moreover, although CAV cooperation control at intersections has undeniable advantages in terms of collision avoidance, constant traffic signal phase and timing (SPAT) control also has advantages in terms of vehicle evacuation, especially for large traffic volumes. In this respect, the comparison between CAV cooperation control and constant traffic SPAT remains deficient. 

Thus, the main contributions of this study are as follows. We propose a two-stage cooperation framework to decouple the longitudinal and lateral control of CAVs. In the first stage, a formation control method is used to guide the CAVs into their target lanes. An iterative framework is developed to solve the multi-vehicle target assignment and path planning problem, preventing the deadlock problem in the single vehicle lane changing algorithm. In the second stage, the CAV arrival time is optimized to increase the traffic efficiency. Specifically, a heuristic graph-based solution is proposed to solve the scheduling problem with low computational burden. As opposed to specially considering lane changing during algorithm design~\cite{lu2019trajectory,hu2019trajectory,xu2020bi}, the decoupling framework has an improved generalization ability in merging with other scheduling methods. In simulations, we compare the proposed algorithm with constant traffic SPAT, and our results verify the effectiveness of the proposed algorithm.

The remainder of this paper is organized as follows. Section~\ref{Sec:ProblemStatement} describes the scenario and introduces the two-stage cooperation framework. Section~\ref{Sec:StageOne} explains the multi-vehicle target assignment and path planning method. Section~\ref{Sec:StageTwo} presents the graph-based CAV scheduling method. The simulation results are discussed in Section~\ref{Sec:Simulation}, and Section~\ref{Sec:Conclusion} concludes the paper.

\section{Two-stage cooperation framework}
\label{Sec:ProblemStatement}
\begin{figure}[!t]
	\centering
	\includegraphics[width=\linewidth]{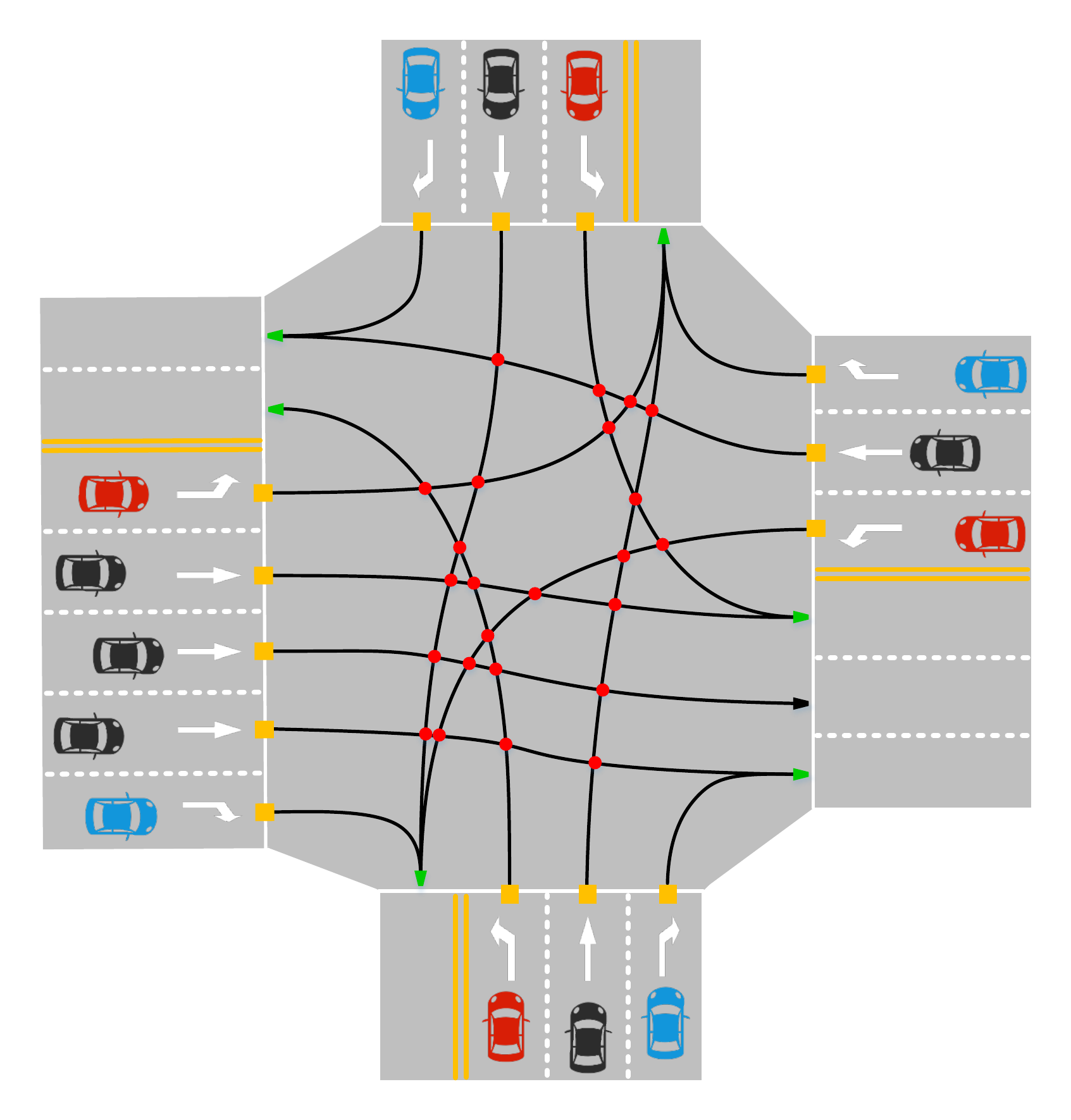}
	\caption{Illustration of the traffic scenario. Vehicles are colored in red, black, and blue to represent turning left, going straight, and turning right, respectively. Red circles, orange squares, and green arrows represent different types of potential collision points.}
	\label{fig:Intersection}
\end{figure}
In~\cite{chen2021conflict}, we introduced a graph-based minimum clique cover (MCC) method applied to an intersection with complex vehicle conflict relationships, as shown in Fig.~\ref{fig:Intersection}. In this previous study, changing lanes was prohibited to guarantee vehicle safety,~\ie, the CAVs are assumed to be in their target lane from the beginning. In this study, the same intersection is used as the traffic scenario; however, lane change is permitted. The CAVs have different destinations; therefore, it is apparent that a scenario where lane changing is permitted is more realistic. Moreover, two assumptions remain to be clarified.

\begin{assumption}
	\label{Ass:Com}
	CAVs transmit their velocities and positions to the central cloud coordinator through ideal wireless communication, ~\eg, V2I communication~\cite{gerla2014internet}, where communication delay and packet loss do not occur. 
\end{assumption}

\begin{assumption}
	\label{Ass:Auto}
	The CAVs are capable of fully autonomous driving, implying that they are assumed to have perfect steering performance in their lane changing and turning behavior.
\end{assumption}

\begin{figure*}[!t]
	\centering
	\includegraphics[width=\linewidth]{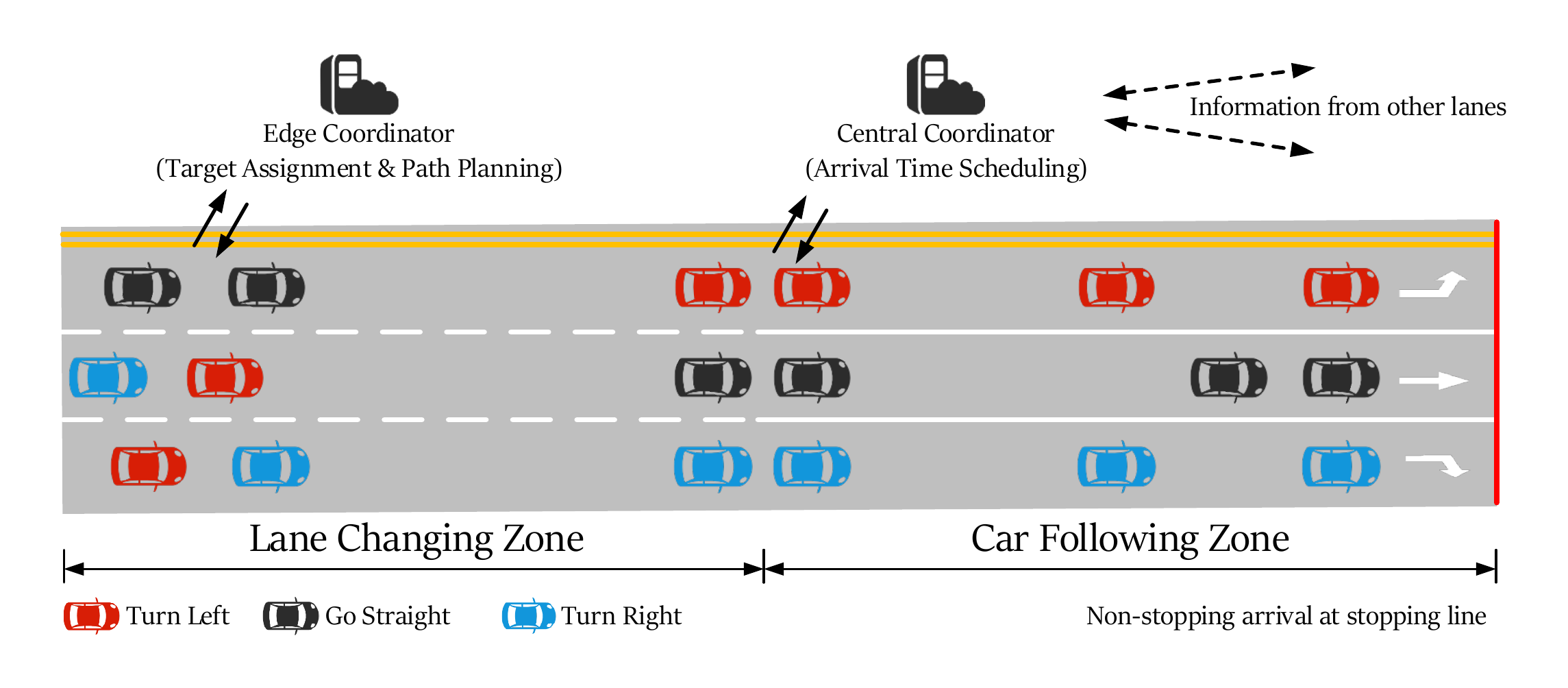}
	\caption{Lane changing and car-following zones. When arriving at the lane changing zone, the CAVs perform lane change maneuvers to arrive at their target lanes. In the car-following zone, the CAVs modify their car-following distances to arrive at the stopping lane according to the scheduled plan.}
	\label{fig:RoadSegmentation}
\end{figure*}

Considering the limit range of V2X communication, previous studies usually set a certain range of the \emph{control zone} to design cooperation algorithms~\cite{chen2020mixed, xu2017v2i, chen2021conflict, xu2018distributed}. In our study, we aim to decouple the longitudinal and lateral control problem in CAVs; therefore, the \emph{control zone} is divided into a \emph{lane changing zone} and a \emph{car-following zone}, as shown in Fig.~\ref{fig:RoadSegmentation}, whose lengths are denoted as $ L_{\mathrm{ctrl}} $, $ L_{\mathrm{LCZ}} $, and $ L_{\mathrm{CFZ}} $, respectively.



Most of the existing studies only considered the trajectory planning of CAVs,~\ie, their longitudinal control, implying that their target lane selection and lane changing behavior was neglected in intersection management. In this study, we propose a two-stage cooperation framework to decouple the longitudinal and lateral control of CAVs, as shown in Fig.~\ref{fig:RoadSegmentation}. The lane changing behavior is strongly related to the surrounding CAVs, and an \emph{edge coordinator} is deployed at the lane changing zone to coordinate the movement of the CAVs to their target lanes. When entering the lane changing zone, task assignment is first accomplished to allocate the preferred lanes to the CAVs. We then develop multi-vehicle path planning to guarantee a collision-free scenario during lane changing. Therefore, when exiting the lane changing zone, all CAVs are aligned to their target positions. As the CAVs enter the car-following zone, the \emph{central coordinator} collects the information from all of them at the intersection and generates an optimal arrival time plan to increase the traffic efficiency. Subsequently, the CAVs modify their car-following distances according to the scheduling result in the car-following zone. Owing to the arrival time of the CAVs being staggered, each CAV travels through the intersection without idling at the stopping line, and thus, vehicle idling is avoided. Notably, in Fig.~\ref{fig:RoadSegmentation}, the lane lines in the lane changing zone are white dashed lines while those in the car-following zone are white solid lines. This implies that lane change is only permitted in the lane changing zone, which corresponds to the current traffic rules.

\begin{remark}
	Decoupling the longitudinal and lateral control of CAVs into two stages has the following advantages. First, many studies have been conducted on the cooperative lane changing behavior of CAVs on straight roads and also arrival time scheduling at intersections. Solving these aforementioned problems separately makes it feasible for us to inherit from these existing studies. Second, although we assume ideal wireless communication in Assumption~\ref{Ass:Com} as in most studies in this field, the effective communication range of the wireless device should not be neglected entirely. Previous experimental tests of the dedicated short-range communications (DSRC) device~\cite{xu2017dsrc} have proved that when the communication distance exceeds $ 500 \mathrm{m} $ and the vehicle velocity is $ 120 \mathrm{km/h}$, the average packet loss rate exceeds $ 20\% $ and the average round-trip time delay exceeds $ 25\mathrm{ms} $. In previous studies, the $ L_{\mathrm{LCZ}} $ and $ L_{\mathrm{CFZ}} $ were usually set to $ 400 - 1000 \mathrm{m}$. Thus, two types of coordinators are deployed to realize intersection management. 
\end{remark}

\section{Multi-vehicle Target Assignment and Path Planning}
\label{Sec:StageOne}

In~\cite{xu2021coordinated}, Cai~\etal proposed a relative coordinate system (RCS) to realize the formation control of CAVs. In~\cite{cai2021formation}, a bi-level conflict-based search (CBS) was further deployed to solve the multi-vehicle collision-free path planning problem. The original methods were applied to a straight road scenario to increase traffic efficiency. We propose that similar methods can be employed to solve the lane changing problem of CAVs.

\subsection{Relative Coordinate System}
\begin{figure*}[!t]
	\centering
	\includegraphics[width=0.8\linewidth]{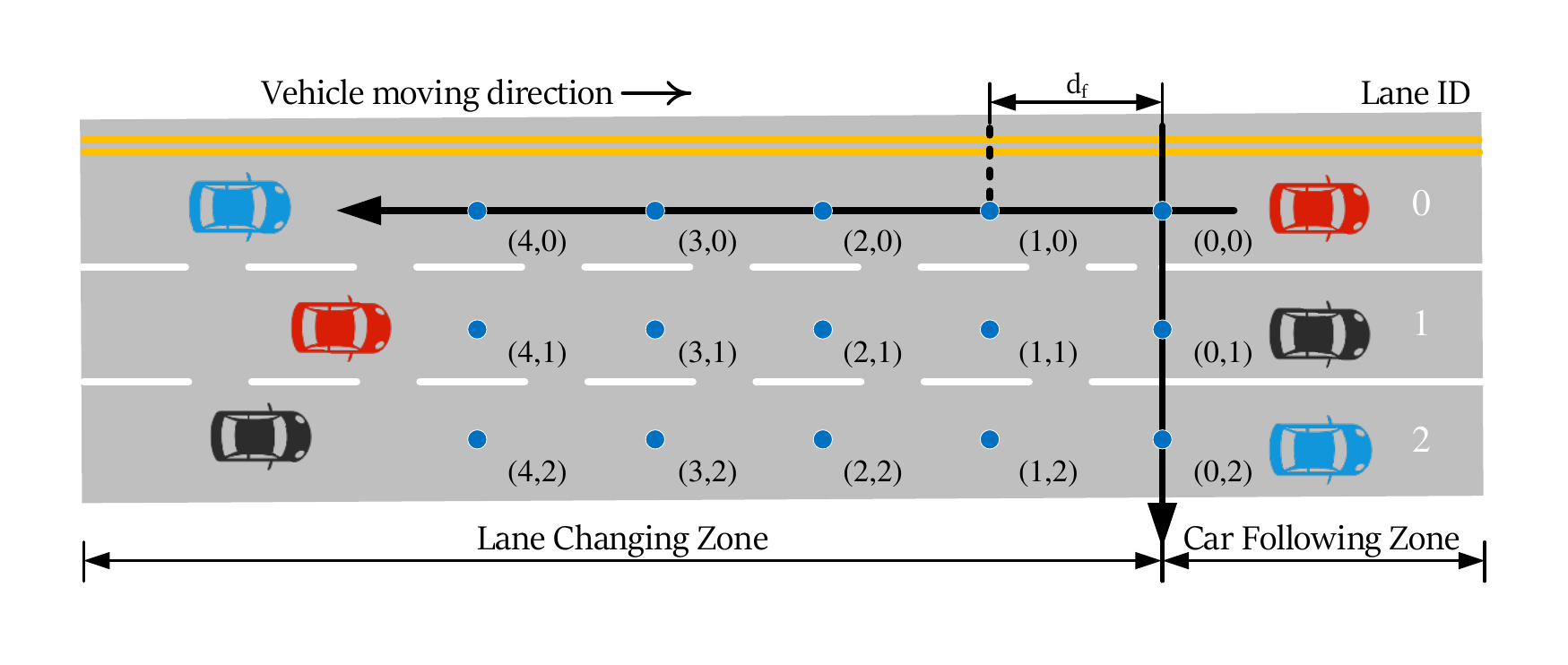}
	\caption{Relative coordinate system (RCS). When entering the lane changing zone, the approaching CAVs are aligned to their nearest coordinates in the RCS. In this zone, the Hungarian algorithm is used to generate the target position assignment and conflict-based search (CBS) is used to find collision-free paths for the CAVs. When exiting the lane changing zone, each CAV adjusts to its target destination.}
	\label{fig:RelativeCoordinateSystem}
\end{figure*}

First, we introduce the RCS in the formation control method. As shown in Fig.~\ref{fig:RelativeCoordinateSystem}, when CAVs arrive at the starting point of the lane changing zone, they are not likely to be in their desired lanes. Hence, the RCS is used to perform the lane change task for the CAVs. Notably, there are two commonly observed formation geometric structures in formation control,~\ie, the parallel and interlaced structures~\cite{kato2002vehicle, marjovi2015distributed}. In~\cite{cai2021formation}, the interlaced structure was selected as the occupation geometric structure to increase the flexibility in multi-vehicle formation coordination. However, in this study, there exists a car-following zone where we adjust the car-following distance. Therefore, we focus on a feasible lane change solution in the lane changing zone. Thus, the parallel structure is selected, where each point in the RCS can be occupied.

When developing conflict-free path planning for a group of CAVs, they are first assigned to the closest points in the RCS. If the two-dimensional position of the CAV $ i $ is $ (x_{i}, \; y_{i}) $, its relative coordinates in the RCS, $ (x^{r}_{i}, \; y^{r}_{i}) $, are determined using
\begin{equation}
	\begin{aligned}
		\min \; &(x_{i} - x^{r}_{i})^2 + (y_{i} - y^{r}_{i})^2,\\
		\mathrm{subject\; to: }\; & x^{r}_{i} \in \mathbb{N}, \; y^{r}_{i} \in \mathbb{N}.
	\end{aligned}
\end{equation}

\subsection{Vehicle-target Assignment}
\label{Sec:Hungary}
After assigning the CAVs to the primary positions in the RCS, their target positions are also assigned based on their destinations. Notably, we use a parallel occupation geometric structure; thus, the CAVs fully occupy the RCS as shown in the car-following zone in Fig.~\ref{fig:RelativeCoordinateSystem}. If there are $ N $ CAVs at primary positions $ (x^{r}_{i}, \; y^{r}_{i}), \; i \leq N, \; i \in \mathbb{N}^+$, there also exist $ N $ target positions $ (x^{t}_{j}, \; y^{t}_{j}), \; j \leq N, \; j \in \mathbb{N}^+ $. Every CAV $ i $ should be assigned to a target position $ j $. 

Before assigning the target positions to the CAVs, the cost of this assignment should be defined. We use the distance from each primary position to each target position as the cost for one CAV, where the Euclidean distance is used as the evaluation index. Hence, we can obtain the cost matrix $ {\mathcal{C}} $ as follows:
\begin{equation}
	\begin{aligned}
		{\mathcal{C}} & = \left[c_{ij}\right] \in \mathbb{R}^{N \times N},\; i,j \in \mathbb{N}^{+},\\
		c_{ij} & = \sqrt{(x^{t}_{k} - x^{r}_{i}) ^{2} + (y^{t}_{k} - y^{r}_{i})},
	\end{aligned}
\end{equation}
where the element in the $ i $-th row and $ j $-th column represents the cost to assign vehicle $ i $ to target $ j $.

Each CAV has a preferred lane; therefore, the preference matrix $ {\mathcal{L}} $ is obtained to define the CAV preferences. 
\begin{equation}
	\begin{array}{c}
		{\mathcal{L}} = \left[l_{ij}\right] \in \mathbb{R}^{N \times N},\; i,j \in \mathbb{N}^{+}, \\
		l_{ij} = \left\{
		\begin{array}{rl}
			1, &\text{if vehicle $ i $ can be assigned to target $ j $,} \\
			M, &\text{otherwise,}
		\end{array}\right.
	\end{array}
\end{equation}
where the element in the $ i $-th row and $ j $-th column represents whether vehicle $ i $ can be assigned to target $ j $ considering their preferred lane. $ M $ is a positive number sufficiently large to prevent vehicle $ i $ from being assigned to target $ j $.

The assignment matrix $ {\mathcal{A}} $, whose element in the $ i $-th row and $ j $-th column represents whether vehicle $ i $ is assigned to target $ j $, is defined as follows:
\begin{equation}
	\begin{array}{c}
		{\mathcal{A}} = \left[a_{ij}\right] \in \mathbb{R}^{N \times N},\; i,j \in \mathbb{N}^{+}, \\
		a_{ij} = \left\{
		\begin{array}{rl}
			1, &\text{if vehicle $ i $ is assigned to target $ j $,} \\
			0, &\text{otherwise.}
		\end{array}\right.
	\end{array}
\end{equation}


Subsequently, the assignment problem can be modeled as a 0-1 integer programming problem:
\begin{equation}
	\label{equ:AssignmentProblem}
	\begin{aligned}
		\min \; & \sum_{i=1}^{N}\sum_{j=1}^{N} \left(c_{ij} \times l_{ij} \times a_{ij} \right), \\
		\mathrm{subject\; to: }\; & \sum_{i=1}^{N} a_{ij} = 1, \\
		& \sum_{j=1}^{N} a_{ij} = 1, \\
		& i,j \in \mathbb{N}^{+},
	\end{aligned}
\end{equation}
where $ N $ is the number of CAVs, $ c_{ij} $ is the cost matrix, $ l_{ij} $ is the preference matrix, and $ a_{ij} $ is the assignment result. 

The Hungarian algorithm~\cite{kuhn1955hungarian} is commonly used to solve assignment problems such as~\eqref{equ:AssignmentProblem}, and we employ this algorithm in our method to generate a feasible assignment with the lowest cost. In the following section, for each group of CAVs, we aim to obtain not only the best assignment but also several sub-optimal assignments. For simplicity, we use $ A_{k} $ to represent the $ k $-th optimal assignment $ {\mathcal{A}} $ and $ C_{k} $ to denote the corresponding cost $ {\mathcal{C}} $.

\subsection{Conflict Types in Lane Changing Behavior}
\begin{figure}[!t]
	\centering
	\subcaptionbox{Node Conflict \uppercase\expandafter{\romannumeral1} \label{fig:Node_Conflict}}
	{\includegraphics[width=0.45\linewidth]{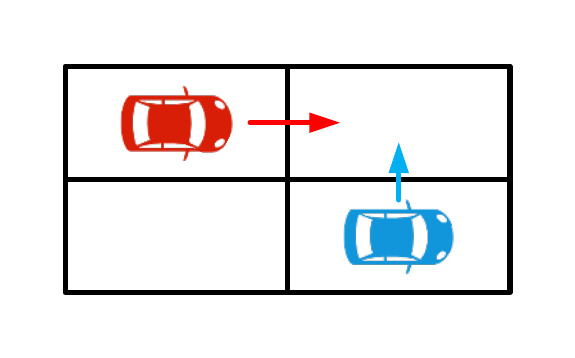}}
	\subcaptionbox{Node Conflict \uppercase\expandafter{\romannumeral2} \label{fig:Node_Conflict2}}
	{\includegraphics[width=0.45\linewidth]{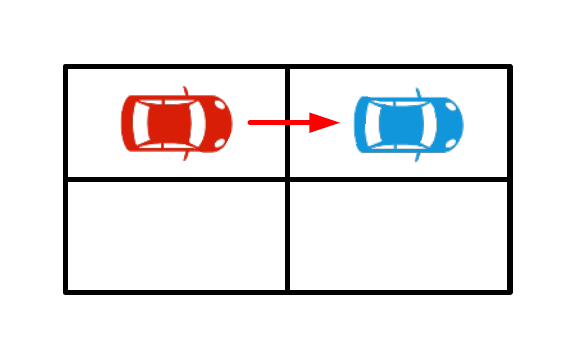}}
	
	\subcaptionbox{Edge Conflict \label{fig:Edge_Conflict}}
	{\includegraphics[width=0.45\linewidth]{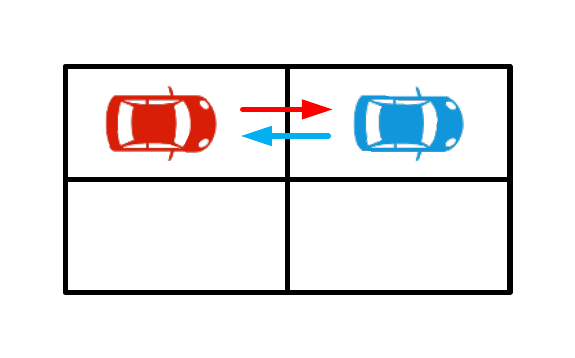}}
	\subcaptionbox{Intermediate Conflict \label{fig:Intermediate_Conflict}}
	{\includegraphics[width=0.45\linewidth]{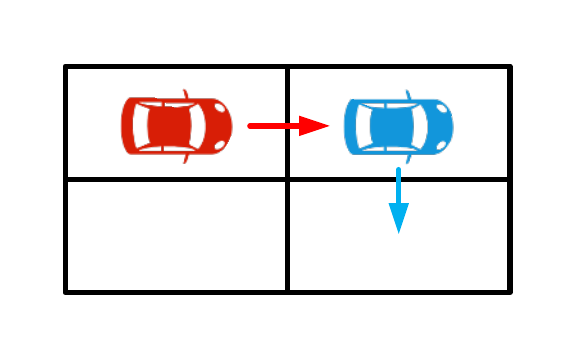}}
	\caption{Conflict types of CAVs in lane change behavior.}
	\label{fig:ConflictRelationships}
\end{figure}

After obtaining the assignment, we need to develop path planning for the CAVs, which should not have conflicts in their trajectories. First, the collision types in the lane changing behavior of CAVs should be clarified. Notably, we regulate the CAVs to move only along orthogonal directions in the RCS,~\ie, movement along diagonal directions in the RCS is prohibited. Movement along diagonal directions in the RCS produces large acceleration and deceleration in lane change behavior, which is hazardous and also increases the computational burden in multi-vehicle path planning. Hence, we constrain the movement of the CAVs to orthogonal directions in the RCS for simplicity.

The conflict types in the path planning of the CAVs are shown in Fig.~\ref{fig:ConflictRelationships}. The first type is the node conflict. As shown in Fig.~\ref{fig:Node_Conflict}, if two CAVs are scheduled into the same coordinates in the RCS, they will inevitably collide. A similar situation occurs in Fig.~\ref{fig:Node_Conflict2}, where one CAV is scheduled into the coordinates of another. The second type is the edge conflict, as shown in Fig.~\ref{fig:Edge_Conflict}, where the scheduled coordinates of two CAVs are interchanged. The third conflict is the intermediate conflict, as shown in Fig.~\ref{fig:Intermediate_Conflict}. In this scenario, one CAV is scheduled to the coordinates of another when the second CAV is on the verge of moving away. Unlike the previous two conflicts, the occurrence of an intermediate conflict depends on the size of the vehicle and the duration of the lane change behavior. For safety concerns, we include this type of conflict in our further analysis.

\subsection{Collision-free Path Planning}

Studies on multi-agent coordination have thoroughly investigated the multi-agent path planning problem. One of the most well-known methods is the CBS ~\cite{sharon2015conflict}, which constructs a constraint tree to obtain the optimal solution in the same manner as single-agent path planning. Other multi-agent path planning methods have also been proposed to solve the problem,~\eg, cooperative path planning~\cite{silver2005cooperative}, the swap method~\cite{luna2011push}, and $ M^{*} $~\cite{wagner2011m}. The primary task of the CAVs in the lane changing zone in this study is to successfully change their lanes to the generated target coordinates. Hence, in this study, the CBS is used to solve the multi-vehicle path planning problem.

$ A^{*} $~\cite{hart1968formal} is a well-known path planning algorithm for a single agent. The input of $ A^{*} $ is the starting point, ending point, and obstacle points. The output is a feasible path from the starting point to the ending point, which avoids the obstacles. Single-agent path planning is not the focus of this study; therefore, we omit the details of the $ A^{*} $ algorithm (hereinafter referred to as A-STAR). Interested readers may refer to~\cite{hart1968formal} for further details.

For each assignment $ A_{k} $ of the $ k $-th iteration, we define $ P_{k} $ as the corresponding collision-free path set. $ P_{k} $ is written as
\begin{equation}
	\begin{aligned}
		P_{k} &= [p_{it}]\in \mathbb{R}^{N \times T},\; i \in \mathbb{N}^{+},t \in \mathbb{N}, \\
		p_{it} &= (x^{r}_{it}, y^{r}_{it}),
	\end{aligned}
\end{equation}
where $ (x^{r}_{it}, y^{r}_{it}) $ are the coordinates of CAV $ i $ in the RCS at time $ t $. Notably, there are $ N $ CAVs in $ A_{i} $; therefore, there are $ N $ rows in $ P_{i} $ representing $ N $ paths. Considering the limited length of the lane changing zone, the maximum number of steps is defined as $ T $. For simplicity, we use $ p_{i} $ to represent the planning path for CAV $ i $ in $ T $ time steps and $ p_{t} $ to represent the coordinates of all $ N $ CAVs at time step $ t $.

The collision-free path set $ P_{k} $ should comply with the following rules. First, the starting points $ (x^{r}_{i0}, y^{r}_{i0}) $ and ending points $ (x^{r}_{iT}, y^{r}_{iT}) $ of CAV $ i $ should correspond to the assignment in $ A_{k} $. Second, at each time step $ t $, the coordinate set $ p_{t} $ should contain no conflict relationships, as shown in Fig.~\ref{fig:ConflictRelationships}. To solve the multi-vehicle path planning problem, the CBS algorithm is proposed in Algorithm~\ref{algo:CBS}.

\begin{algorithm}[t]
	\caption{Conflict-based Search}
	\label{algo:CBS} 
	\begin{algorithmic}[1] 
		\Require{The $ k $-th assignment $ A_{k} $ for a group of $ N $ CAVs} 
		\Ensure{Collision-free planning path set $ P_{k} $}
		
		\State Initialize conflict point set $ F_{i} = \varnothing, \, i \in \mathbb{N}^{+} $
		
		\While{True}
			\For{each CAV $ i $ in assignment $ A_{k} $} 
				\State $ p_{i} $ = A-STAR($ (x^{r}_{i0}, y^{r}_{i0}) $, $ (x^{r}_{iT}, y^{r}_{iT}) $, $ F_{i} $)
				\State update the path $ p_{i} $ in $ P_{k} $
			\EndFor
		
			\If{there exists conflict $ f $ for CAV $ i \in \mathbb{N}^{+}$ in $ P_{k} $}
				\State $ F_{i} = F_{i} \cup f $
			\Else
				\State \Return $ P_{k} $
			\EndIf
		\EndWhile
	\end{algorithmic} 
\end{algorithm}

The essential part of the CBS algorithm is examining the conflicts $ f $ in the current planning path set $ P_{k} $, adding them into the corresponding conflict point set $ F_{i} $ for CAV $ i $, and executing the next round of path planning. The iteration ends when there exist no conflicts in $ P_{k} $,~\ie, the planning path set $ P_{k} $ is collision-free. It can be argued that this method may not find the collision-free path planning solution for a large number of CAVs. However, the feasibility of the algorithm also depends on the maximum number of steps. In the worst case, provided that the maximum time step $ T $ is sufficiently large, the CAVs can always find a collision-free path planning solution. Notably, because the assignment $ A_{k} $ does not consider the conflict while the planning path $ P_{k} $ is collision-free, the cost of assignment $ A_{k} $ $ C_{k} $ never exceeds that of the planning path $ P_{k} $, which is denoted as $ C'_{k} $. We interpret that this non-decreasing property of the CBS guarantees the optimality of the algorithm.

\subsection{Iterative Solution for Assignment and Path Planning}
\label{Sec:PathPlanning}
\begin{figure*}[!t]
	\centering
	\subcaptionbox{Algorithm interpretation\label{fig:CBS_algorithm}}
	{\includegraphics[width=0.52\linewidth]{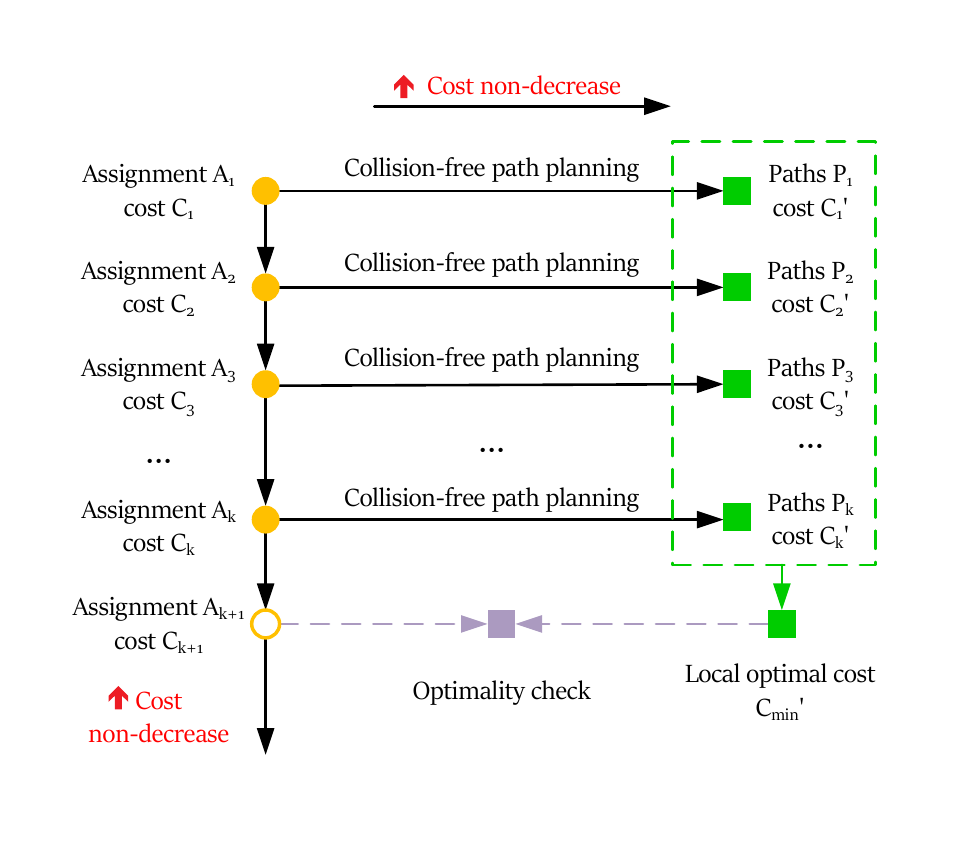}}
	\subcaptionbox{Algorithm flow chart \label{fig:CBS_flow}}
	{\includegraphics[width=0.45\linewidth]{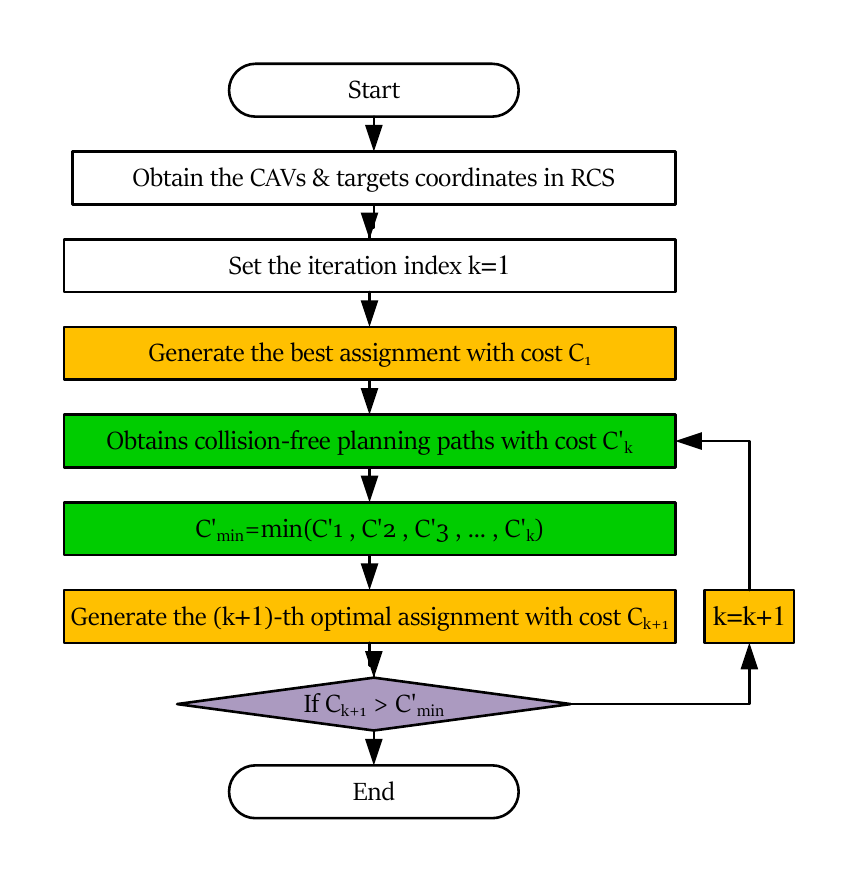}}
	\caption{Iterative framework for solving the target assignment and path planning problem. Yellow blocks represent target assignment sections. Green blocks stand for path planning sections. Purple ones are the optimality check process.}
	\label{fig:ConflictBasedSearching}
\end{figure*}

In the previous sections, we have shown how we generate the best assignment $ A_{1} $ according to the initial positions and target lanes of the CAVs. We have also introduced how we obtain a collision-free planning path set $ P_{1} $ from the best optimal assignment $ A_{1} $. However, $ P_{1} $ is the local optimal solution generated from $ A_{1} $. We have not checked all the assignments; therefore, a global optimal solution has not been found. Thus, we propose an iterative framework to find this global optimal solution, as shown in Fig.~\ref{fig:ConflictBasedSearching}. In the figure, we use the colors yellow to illustrate the part of the target assignment solved using the Hungarian algorithm (Section~\ref{Sec:Hungary}) and green to denote the part of the multi-vehicle path planning solved using A*-based CBS (Section~\ref{Sec:PathPlanning}).

As shown in Fig.~\ref{fig:CBS_flow}, the CAVs are first assigned to their target coordinates. The Hungarian algorithm provides the optimal assignment $ A_{1} $ with the lowest cost $ C_{1} $,~\ie, each CAV is assigned to unique target coordinates based on their permitted lanes with the shortest total distance. We then develop the initial path planning for each CAV and obtain the initial collision-free path $ P_{1} $ with the cost $ C_{1}' $. However, we cannot claim that $ P_{1} $ is the global optimal solution because $ P_{k} $, generated using another $ A_{k} $, could have a lower cost. Thus, we generate the second-best assignment $ A_{2} $ using the Hungarian algorithm, which has the lowest cost except for $ A_{1} $. $ A_{2} $ is not the best assignment; therefore, the cost $ C_{1} \geq C_{2} \geq C_{k}, \, k>2 $. It is apparent that if $ C_{2} > C_{1}' $, $ P_{1} $ is the global optimal solution because the path $ C_{2}' $ generated from $ P_{2} $ has a larger cost than $ C_{2} $ and all the remaining assignments also have larger costs. In other circumstances, if $ C_{2} \leq C_{1}' $, we further generate the planning path $ P_{2} $ and the third-best optimal assignment $ A_{3} $, and the iteration continues. The determination of the optimality is colored in purple in Fig.~\ref{fig:ConflictBasedSearching}.

\begin{theorem}
	The planning path set $ P_{\min} $ generated in Fig.~\ref{fig:CBS_flow} is the global optimal solution,~\ie, the cost $ C_{\min}' $ is the lowest cost of the collision-free path set. 
\end{theorem}

\begin{IEEEproof}
	Without loss of generality, we set $ C_{i}' = \min(C_{1}', C_{2}', C_{3}', \dots, C_{k}') $, implying that the generated planning path set $ P_{\min} = P_{i} $. If there exists a planning path set $ P_{j} $ with a lower cost $ C_{j}' \leq C_{i}' $, it must be generated from the assignment $ A_{j} $ and $ C_{j}'\geq C_{j} $. $ P_{j} $ is not explored in the previous path plannings; therefore, it must be generated from an assignment that has not been explored,~\ie, $ j > k $. Thus, we have $ C_{j} \geq C_{k} $. From the algorithm, we have $ C_{k} > C_{\min}' = C_{i}' $; therefore, we summarize the inequality equations as $ C_{j}'\geq C_{j} \geq C_{k} \geq C_{i}' $, which contradicts the assumption $ C_{j}' \leq C_{i}' $. Therefore, $ C_{i}' $ is the lowest cost and $ P_{\min} $ is the global optimal solution.
\end{IEEEproof}

\subsection{Vehicle Control}
After obtaining the planning path set $ P_{\min} $, generated from the iterative solution method, each CAV has a collision-free path to its target position. CAV control is not the focus in this study; therefore, we simplify the CAV control process to ease the calculation burden. In Assumption~\ref{Ass:Auto}, we have assumed that the CAV has perfect steering performance, which exempts us from a lateral controller design. Therefore, given adequate time, CAVs are able to execute the planned movement in one step. Hence, a second-order vehicle model is used, and the deviation from the current position of CAV $ i $ to its target position at time $ t $ is defined as follows:
\begin{equation}
	\begin{array}{l}
		\delta_{\mathrm{p}}^{(i,t)} = x_{i,t} - x^{r}_{it},\\
		\delta_{\mathrm{v}}^{(i,t)} = v_{i,t}(t)-v_{\mathrm{p}},
	\end{array}
\end{equation}
where $ (x_{i,t}, y_{i,t}) $, and $ v_{i,t}(t) $ represent the position and velocity, respectively, of CAV $ i $ at time $ t $, and $ v_{\mathrm{p}} $ is the designed platoon velocity. A linear feedback controller is designed as follows:
\begin{equation}
	\label{equ:Feedback_1}
		u_{i}=-k_{p} \delta_{\mathrm{p}}^{(i, j)}-k_{v} \delta_{\mathrm{v}}^{(i, j)}.
\end{equation}

The design of the RCS coordinates and the controller also considers the CAV control in the car-following zone, which is described in Section~\ref{Sec:VehicleControl}. We obtain the lane change timing for each CAV from the planning path set $ P_{\min} $; however, the lane changing behavior is accomplished by the lane changing model in the traffic simulator SUMO~\cite{erdmann2015sumo}. Other key control parameters are listed in Table~\ref{tab:Parameters}.

%

\section{CAV Scheduling at Unsignalized Intersections}
\label{Sec:StageTwo}
As shown in Fig.~\ref{fig:RoadSegmentation}, the CAVs run in their target lanes after driving through the lane changing zone. The trajectories of the CAVs intersect in the middle of the intersection; thus, in the car-following zone, the CAVs have to schedule collision-free arrival plans. In~\cite{xu2018distributed}, a method using a \emph{virtual platoon} was proposed, which projects the CAVs from different lanes onto a virtual lane. Thus, CAVs from different lanes can drive through the intersection as if they were in the same lane. Therefore, the central coordinator only needs to schedule the CAVs and transmit their target positions in the virtual platoon. Subsequently, CAVs can be controlled in terms of platooning behavior.

\subsection{Vehicle Model}
\begin{figure}[!t]
	\centering
	\includegraphics[width=\linewidth]{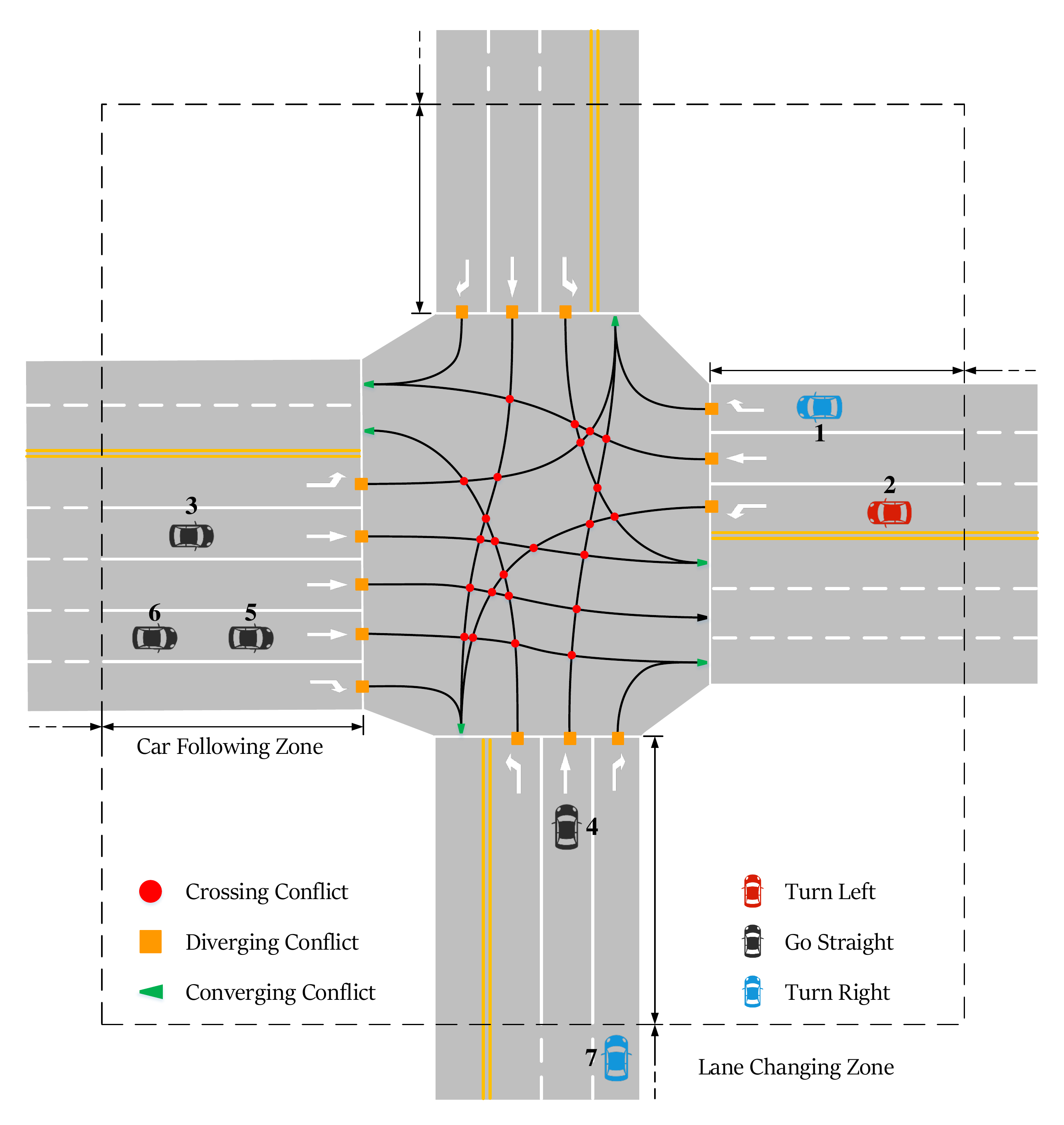}
	\caption{Traffic scenario in the car-following zone. Potential conflict points are denoted by red circles, orange squares, and green triangles. CAVs are colored in red, black, or blue based on their different destinations.}
	\label{fig:MCC_Scenario}
\end{figure}

The incoming CAVs are indexed from $ 1 $ to $ N $ according to their arrival sequence in the car-following zone, as shown in Fig.~\ref{fig:MCC_Scenario}. For each vehicle $ i \left(i \leq N, i \in \mathbb{N}^+\right) $, the second-order dynamic model is given by

\begin{equation}
	\label{equ:StateSpaceEquation}
	\left\{
	\begin{array}{l}
		\dot{x_{i}}(t) = v_{i}(t), \\
		\dot{v_{i}}(t) = u_{i}(t).
	\end{array} \right.
\end{equation}
$ u_{i}(t) $ represents the input of vehicle $ i $ at time $ t $. There also exist velocity and acceleration constraints on the vehicle.
\begin{equation}
	\label{equ:Constraints}
	\begin{aligned}
		0 &\le v_{i} \le v_{\max}, \\
		u_{\min} &\le u_{i} \le u_{\max}.
	\end{aligned}
\end{equation}

\subsection{Conflict Analysis}

The CAVs from different lanes form multiple conflict points. Despite the complicated conflict scenarios, they can be classified into the following conflict modes. Without loss of generality, several CAVs are selected in Fig.~\ref{fig:MCC_Scenario} to illustrate the conflict relationship. We define four conflict types in this study.

\begin{enumerate}
	\item \emph{Crossing Conflict}: CAVs from different lanes have the potential to collide while crossing the conflict points, indicated by the 24 red circles. For example, CAV $ 2 $ and CAV $ 3 $ have a crossing conflict point.
	\item \emph{Diverging Conflict}: The CAVs change lanes in the lane changing zone, and lane changing and passing are not permitted in the car-following zone. Thus, vehicles on the same lane cannot pass the intersection simultaneously, as indicated by the 14 orange squares. For example, CAV $ 5 $ and CAV $ 6 $ have a diverging conflict point.
	\item \emph{Converging Conflict}: Vehicles from different lanes cannot drive into the same lane simultaneously, as shown by the six green arrows. For example, CAV $ 1 $ and CAV $ 4 $ have a converging conflict point.
	\item \emph{Reachability Conflict}: CAVs cannot pass the intersection simultaneously because of the acceleration and velocity constraints, regardless of whether or not they have the abovementioned conflicts. CAV $ 1 $ and CAV $ 7 $ have a reachability conflict.
\end{enumerate}

The first three conflict types are route conflicts~\cite{roess2004traffic}, where CAVs have intersections in their trajectories along their paths. The fourth conflict type is caused by the velocity and acceleration constraints on the CAVs. For instance, in Fig.~\ref{fig:MCC_Scenario}, CAV $ 7 $ arrives at the car-following zone when CAV $ 1 $ nearly reaches the stop line with the designed virtual platoon velocity. In this case, CAV $ 7 $ cannot catch up with CAV $ 1 $ at the stop line, regardless of whether or not they have any conflict relationships. This is because of the constraints on the vehicle velocity and acceleration.

Expanding equation~\eqref{equ:StateSpaceEquation} and~\eqref{equ:Constraints}, we obtain the evaluation condition as follows:
\begin{equation}
	\label{equ:Reachability}
	\frac{L_\mathrm{prec}}{v_\mathrm{p}} < \frac{L_\mathrm{CFZ}}{v_{\max}} + \frac{v_{\max}}{2u_{\max}},
\end{equation}
where $ L_\mathrm{prec} $ and $ L_\mathrm{CFZ} $ are the distances from the stop line to the preceding CAV and to the beginning of the car-following zone, respectively, and $ v_\mathrm{p} $ is the designed virtual platoon velocity. Namely, the CAVs that are very close to the intersection should not be considered in the scheduling of new incoming CAVs. Notably, most of the existing studies assumed that the CAVs can reach the stop line under all circumstances,~\ie, the reachability conflict was ignored.

We define different conflict sets to describe the conflict relationship of the CAVs. For each CAV $ i \left(i \leq N, i \in \mathbb{N}^+\right) $, the crossing, diverging, converging, and reachability sets are defined as $ \mathcal{C}_{i} $, $ \mathcal{D}_{i} $, $ \mathcal{V}_{i} $, and $ \mathcal{R}_{i} $, respectively. Notably, because the conflict sets are determined when the CAV reaches the car-following zone, the CAV indexes in the conflict sets are smaller than those of the CAVs at the border of the car-following zone,~\ie, the elements in the conflict sets satisfy~\eqref{equ:ElementRelationship}.
\begin{equation}
	\label{equ:ElementRelationship}
	i < j, \, \mathrm{if}\ i \in \mathcal{C}_{j} \cup \mathcal{D}_{j} \cup \mathcal{V}_{j} \cup \mathcal{R}_{j}.
\end{equation}

\subsection{Description of Graph-based Conflicts}
\label{Sec:ConflictDescription}
Based on the conflict set analysis, we further define a conflict directed graph (CDG) $ \mathcal{G}_{N+1} $ to represent the conflict relationship between the CAVs.

\begin{definition}[Conflict Directed Graph]
	\label{def:CDG}
	The CDG is denoted as $ \mathcal{G}_{N} = \left( \mathcal{V}_{N},\mathcal{E}_{N} \right) $. If there are $ N $ CAVs in the car-following zone, we have the node set $ \mathcal{V}_{N} = \{1,2,\dots,N\} $. The unidirectional edge set is defined as $ \mathcal{E}_{N}^{u}=\{(i,j) \mid i \in \mathcal{D}_{j} \cup \mathcal{R}_{j}\} $, and the bidirectional edge set is defined as $ \mathcal{E}_{N}^{b}=\{(i,j) \mid i \in \mathcal{C}_{j} \cup \mathcal{V}_{j}\} $. The edge set is the union of these two sets as $ \mathcal{E}_{N} = \mathcal{E}_{N}^{u} \cup \mathcal{E}_{N}^{b} $.
\end{definition}

\begin{figure}[!t]
	\centering
	\subcaptionbox{Conflict Directed Graph\label{fig:ConflictDirectedGraph}}
	{\includegraphics[width=0.53\linewidth]{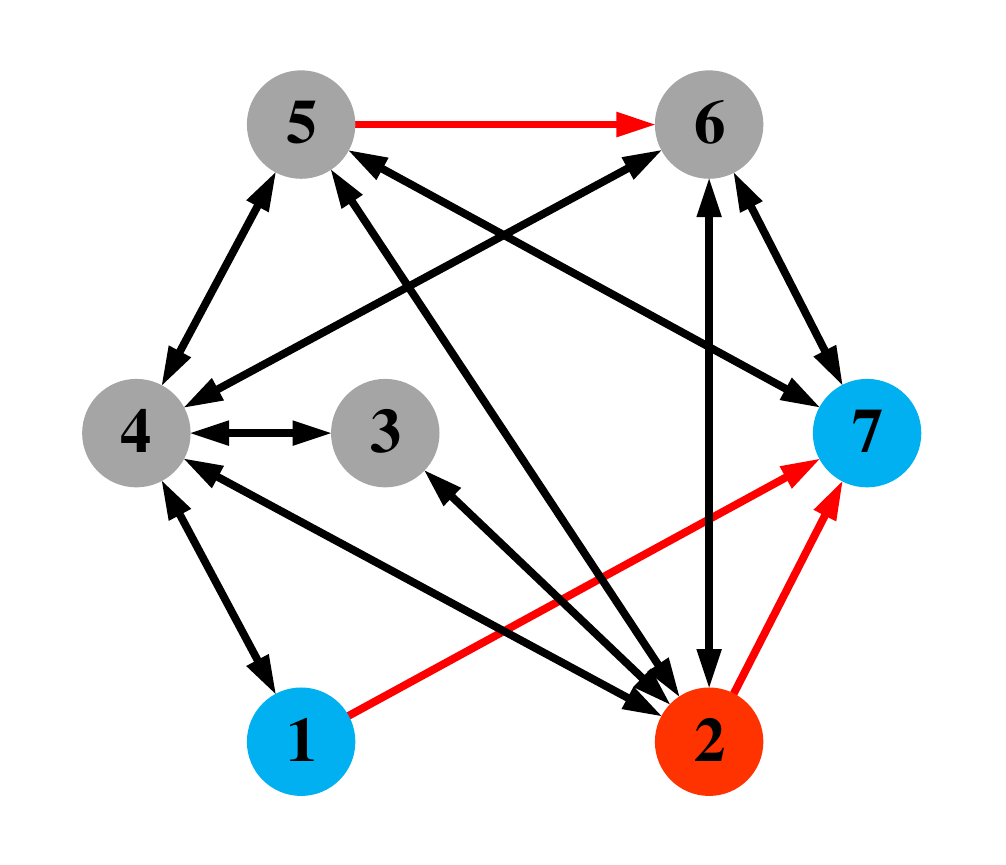}}
	\subcaptionbox{Coexisting Undirected Graph\label{fig:CoexistUndirectedGraph}}
	{\includegraphics[width=0.41\linewidth]{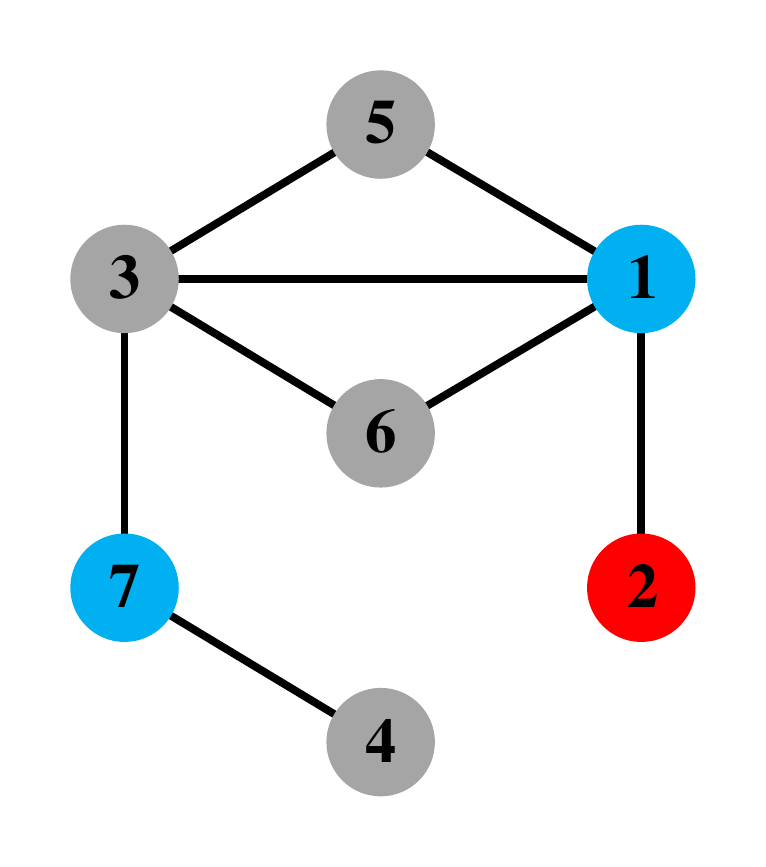}}
	\caption{Fig.~\ref{fig:ConflictDirectedGraph} is the conflict directed graph (CDG). The red unidirectional edges represent the diverging and reachability conflicts, whereas the black bidirectional edges represent the crossing and converging conflicts. Fig~\ref{fig:CoexistUndirectedGraph} is the coexisting undirected graph (CUG), which is the complement graph of the CDG and describes the coexistence relationship of the vehicles.}
	\label{fig:ConflictDirectedScenario}
\end{figure}

The CDG of the scenario in Fig.~\ref{fig:MCC_Scenario} is depicted in Fig.~\ref{fig:ConflictDirectedGraph}. The nodes in the CDG represent the CAVs in the car-following zone. The red unidirectional edges represent the diverging and reachability conflicts. The existence of a unidirectional edge $ (i,j) $ implies that CAV $ j $ is not permitted to pass CAV $ i $ or CAV $ j $ is unable to catch up to CAV $ i $ because it satisfies~\eqref{equ:Reachability}. Thus, CAV $ j $ cannot reach the intersection earlier than CAV $ i $. The black bidirectional edges denote the crossing and converging conflicts, implying that the arrival sequences of CAVs $ i $ and $ j $ can be interchanged. 

It is straightforward that the CDG describes all the conflict relationships of the CAVs. From a different perspective, another method to describe the conflict relationships of the CAV is to describe their coexistence relationships.

\begin{definition}[Coexisting Undirected Graph]
	\label{def:CUG}
	The coexisting undirected graph (CUG) is defined as the complement graph of the CDG $ \mathcal{G}_{N} $. Thus, $ \overline{\mathcal{G}}_{N} = \left(\overline{\mathcal{V}}_{N},\overline{\mathcal{E}}_{N}\right) $, where $ \overline{\mathcal{V}}_{N} = \mathcal{V}_{N} $, $ \overline{\mathcal{E}}_{N} = \{(i,j) \mid i,j \in \overline{\mathcal{V}}_{N} , i \neq j, \mathrm{and} \ (i,j) \notin \mathcal{E}_{N}\} $.
\end{definition}

In this scenario, the CDG is depicted in Fig.~\ref{fig:ConflictDirectedGraph} and the CUG in Fig.~\ref{fig:CoexistUndirectedGraph}. The CDG edge $ \mathcal{E}_{N} $ implies that two CAVs have conflicts and the CUG is the complement graph of the CDG; therefore, the CUG edge $ \overline{\mathcal{V}}_{N} $ implies that the two CAVs are conflict-free,~\ie, they can pass through the intersection simultaneously.


\subsection{Minimum Clique Cover Scheduling}
Generally, all the scheduling methods at unsignalized intersections are designed to find a high-efficiency collision-free passing order of the CAVs,~\ie, to schedule the arrival sequence based on their conflict relationships described in Section~\ref{Sec:ConflictDescription}. In this study, we focus on scheduling the CAVs based on the CUG defined in Definition~\ref{def:CUG}. The CUG describes the coexistence of the CAVs, which may pass the intersection simultaneously. In graph theory, a clique is suitable for describing the coexistence relationship of CAVs. The definition of the clique is shown in Definition~\ref{def:clique}.

\begin{definition}[Clique\cite{luce1949method}]
	\label{def:clique}
	A clique $ C $ in an undirected graph $ G = (V, E) $ is a subset of the nodes, $ C \subseteq V, $ such that every two distinct nodes are adjacent. This is equivalent to the condition that the subgraph of $ G $ induced by $ C $ is a complete graph.
\end{definition}

Considering the cliques in the CUG $ \overline{\mathcal{G}}_{N} $, the CAVs in one clique are conflict-free,~\ie, they can pass through the intersection simultaneously. Therefore, the objective is to find the minimum number of vehicle groups in the CUG,~\ie, the minimum number of cliques covering all the nodes in the CUG $ \overline{\mathcal{G}}_{N} $. Therefore, we define the MCC problem as follows.

\begin{definition}[Minimum Clique Cover (MCC)~\cite{karp1972reducibility}]
	\label{def:MCC}
	A clique cover of a graph $ G=(V,E) $ is a partition of $ V $ into $ k $ disjoint subsets $ V_{1},V_{2},\dots, V_{k}$ such that for $ 1 \leq i \leq k $, the subgraph induced by $ V_{i} $ is a clique,~\ie, a complete graph. The MCC number of $ G $ is the minimum number of subsets in a clique cover of $ G $, denoted as $ \theta(G) $.
\end{definition}

The MCC number $ \theta(\overline{\mathcal{G}}_{N}) $ of the CUG represents the minimum number of cliques covering it. The cliques in the CUG represent the CAVs that can pass through the intersection simultaneously; therefore, these CAVs in the same clique can be scheduled to simultaneously drive through the intersection. For example, considering the CUG in Fig.~\ref{fig:CoexistUndirectedGraph}, the MCC number $ \theta(\overline{\mathcal{G}}_{N}) = 4 $, and the corresponding cliques are listed in Table~\ref{tab:MCC}. We conclude that for an arbitrary intersection scenario, the coexistence relationship of the incoming CAVs is depicted in the CUG $ \overline{\mathcal{G}}_{N} $. Thus, the MCC number $ \theta(\overline{\mathcal{G}}_{N}) $ represents the possible minimum passing order solution.

\begin{table}[!t]
	\renewcommand{\arraystretch}{1.3}
	\centering
	\begin{tabular}{c|c|c|c|c}
		\hline
		\diagbox[width=8em]{\textbf{Solutions}}{\textbf{Subsets}} & $ {V}_{1} $ & $ {V}_{2} $  & $ {V}_{3} $ & $ {V}_{4} $ \\
		\hline
		\bfseries 1 & $ \{1,3,5\} $ & $ \{4,7\} $ & $ \{2\} $ & $ \{6\} $ \\
		\bfseries 2 & $ \{1,3,6\} $ & $ \{4,7\} $   & $ \{2\} $ & $ \{5\} $\\
		\bfseries 3 & $ \{1,2\} $ & $ \{3,5\} $ & $ \{4,7\} $ & $ \{6\} $\\
		\bfseries 4 & $ \{1,2\} $ & $ \{3,6\} $ & $ \{4,7\} $ & $ \{5\} $\\
		\bfseries 5 & $ \{1,5\} $ & $ \{3,6\} $ & $ \{4,7\} $ & $ \{2\} $\\
		\bfseries 6 & $ \{1,6\} $ & $ \{3,5\} $ & $ \{4,7\} $ & $ \{2\} $\\
		\hline
	\end{tabular}
	\caption{Possible MCC solutions of Fig.~\ref{fig:CoexistUndirectedGraph}. Note that $ \theta(\overline{\mathcal{G}}_{N}) = 4$ in this graph.}
	\label{tab:MCC}
\end{table}

\begin{figure}[!t]
	\centering
	\includegraphics[width=0.7\linewidth]{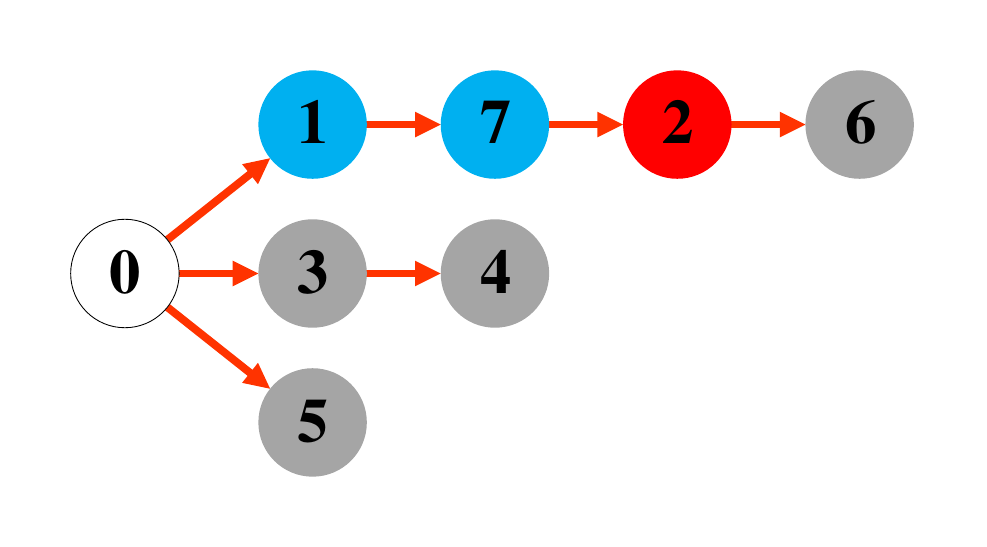}
	\caption{Spanning tree generated by the MCC method. The CAVs in the same layer are collision-free and can pass the intersection simultaneously. Therefore, the overall depth of the spanning tree is $ 4 $, implying that all the CAVs can pass the intersection in the platoon passing time of $ 4 $ vehicles.}
	\label{fig:Spanning_Tree}
\end{figure}

If we consider solution $ 1 $ from Table~\ref{tab:MCC} and place CAVs of the same clique into the same layer, a spanning tree $ \mathcal{G}_{N+1}' $ is generated as shown in Fig~\ref{fig:Spanning_Tree}. Node $ 0 $ is the virtual leading vehicle and the CAVs in the same layer are arranged to simultaneously pass the intersection. Notably, $ \theta(\overline{\mathcal{G}}_{N}) = 4 $ in the CUG in Fig.~\ref{fig:CoexistUndirectedGraph}, implying that the spanning trees generated by the MCC solutions have a minimum layer of $ 4 $. This also indicates that the theoretical evacuation times of these solutions are the same,~\ie, the theoretical values of $ t_\mathrm{evac} $ are the same. Thus, the evaluation index of $ t_\mathrm{evac} $ in~\eqref{equ:EvacuationTime} corresponds to $ \theta(\overline{\mathcal{G}}_{N}) $, which is the global optimal passing order considering the evacuation time. In addition, we further consider the average travel time delay (ATTD) among these solutions as a secondary index. The definition of $ t_\mathrm{ATTD} $ in~\eqref{equ:TravelTimeDelay} can be rewritten in graphical terms as
\begin{equation}
	\label{equ:ATTD_MCC}
	\begin{aligned}
		\min \; & \sum_{i=1}^{k} d_{i}|V_{i}|, \, i \in \mathbb{N}^+,\\
		\mathrm{subject\; to: }\; & \sum_{i=1}^{k} V_{i} = N,
	\end{aligned}
\end{equation}
where $ V_{i} $ represents the nodes in CUG $ \overline{\mathcal{G}}_{N} $ and $ d_{i} $ is the spanning tree depth of the each CAV node $ V_{i} $.

In this scenario, it is evident that the subsets $ V_{i} $ should be arranged in the descending order to decrease the average $ d_{i} $ of $ N $ CAVs in~\eqref{equ:ATTD_MCC}. For example, in Table~\ref{tab:MCC}, we prefer to choose solution $ 1 $, and the corresponding optimized spanning tree is scheduled as $ \{1,3,5\} \rightarrow \{4,7\} \rightarrow \{2\} \rightarrow \{6\}$, as shown in Fig.~\ref{fig:Spanning_Tree}.

The CUG $ \overline{\mathcal{G}}_{N} $ only contains the coexistence information of the CAVs, whereas the passing order of the CAVs in the same lane should be strictly according to their relative positions. If the MCC method generates solutions that cannot be directly executed,~\eg, solution $ 2 $ in Table~\ref{tab:MCC}, the spanning tree is generated as $ \{1,3,6\} \rightarrow \{4,7\} \rightarrow \{2\} \rightarrow \{5\} $. However, this solution is not feasible because CAV $ 5 $ is ahead of CAV $ 6 $. We have proved that this can be solved by exchanging the unfeasible sequence~\cite{chen2021conflict}, as shown in Lemma~\ref{Pro:1}.

\begin{lemma}
	\label{Pro:1}
	If the MCC method leads to an unfeasible solution, where CAVs $ i $ and $ j $ have conflicting trajectories, exchanging the positions of $ i $ and $ j $ solves the conflict. The new solution is also an MCC solution.
\end{lemma}

The MCC problem is proved to be an NP-hard problem, which is difficult to solve during real deployment, especially when the CAV number $ N $ is large. Therefore, we apply a practical approach to solve the problem heuristically, as shown in lines~\ref{algo:MCCBegin} to~\ref{algo:MCCEnd} of Algorithm~\ref{algo:MCC}. The MCC problem of $ G $ is proved to be reduced to the graph coloring problem of $ \overline{G} $~\cite{garey1979guide}, and there are numerous heuristic methods to solve the graph coloring problem. As mentioned earlier, we intend to find the solutions with larger cliques. Thus, we first generate a node sequence $ K=(v_{1},v_{2},\dots,v_{N}) $ using a breadth-first search (BFS). Then, we greedily assign node $ v_{i} $ the smallest possible color,~\ie, the clique index according to the node sequence $ K $, forming the subset cliques $ V_{1},\dots,V_{k}$. Line~\ref{algo:MCC:Spanning} depicts the spanning process, which arranges the cliques into a spanning tree and has a constant calculation time.

\begin{algorithm}[tb]
	\caption{Minimum Clique Cover Method}
	\label{algo:MCC} 
	\begin{algorithmic}[1] 
		\Require{Coexisting Undirected Graph $\overline{\mathcal{G}}_{N} = \left(\overline{\mathcal{V}}_{N},\overline{\mathcal{E}}_{N}\right) $} 
		\Ensure{Spanning Tree $ \mathcal{G}_{N+1}' = \left( \mathcal{V}_{N+1},\mathcal{E}_{N+1}' \right) $}
		\State Calculate the complement graph $\overline{\overline{\mathcal{G}}}_{N} = \mathcal{G}_{N}$ \label{algo:MCCBegin}
		\State Find the breadth-first search sequence $ K=(v_{1},v_{2},\dots,v_{N}) $ of the nodes in $ \mathcal{G}_{N} $
		\For{each node $ v_{i} $ of $ \mathcal{G}_{N} $ in the sequence $ K $} 
		\State assign node $ v_{i} $ the smallest possible clique index
		\EndFor \label{algo:MCCEnd}
		\State Rank $ V_{1}, V_{2}, \dots,V_{k}$ in the descending order and obtain the spanning $ \mathcal{G}_{N+1}' $ \label{algo:MCC:Spanning} 
		\State Exchange the conflicting CAVs of $ \mathcal{G}_{N+1}' $ in the same lane if necessary
		\State \Return $ \mathcal{G}_{N+1}' $
	\end{algorithmic} 
\end{algorithm}

\subsection{Vehicle Control}
\label{Sec:VehicleControl}
A virtual platoon is constructed from the spanning tree in Fig.~\ref{fig:Spanning_Tree}; hence, the controller of CAV $ i $ is related to the communication topology design. Owing to the page limit, we only present the basic controller design here. Notions $ a_{ij} $, $ q_{ij} $, and $ l_{ij} $ are related to the designed predecessor--leader following (PLF) topology. Interested readers may refer to~\cite{chen2021conflict} for further details.

First, a union set $ \mathbb{I}_{i} $ is defined to describe the information exchange of CAV $ i $ as follows:
\begin{equation}
	\mathbb{I}_{i}=\left\{j \mid a_{i j}=1\right\} \cup \left\{0 \mid q_{i i}=1\right\}.
\end{equation}

The distance and velocity errors are defined as
\begin{equation}
	\begin{array}{l}
		\delta_{\mathrm{p}}^{(i,j)} = p_{j}(t)-p_{i}(t)-d_\mathrm{f}\left(d_{j}-d_{i}\right) \\
		\delta_{\mathrm{v}}^{(i,j)} = v_{i}(t)-v_{j}(t)
	\end{array}, \, j \in \mathbb{I}_{i},
\end{equation}
where $ d_\mathrm{f} $ is the car-following distance, $ \delta_{\mathrm{p}}^{(i,j)} $ is the car-following distance error of CAV $ i $, and $ \delta_{\mathrm{v}}^{(i,j)} $ is the car-following velocity error considering all the CAVs in $ \mathbb{I}_{i} $.

A linear feedback controller is designed as follows:
\begin{equation}
	\label{equ:Feedback_2}
	\begin{aligned}
		u_{i} 
		&=-\sum_{j \in \mathbb{I}_{i}} k_{p} \delta_{\mathrm{p}}^{(i, j)}-\sum_{j \in \mathbb{I}_{i}} k_{v} \delta_{\mathrm{v}}^{(i, j)},
	\end{aligned}
\end{equation}
where $ k_{p} $ and $ k_{v} $ are the feedback gains of the distance and velocity errors of CAV $ i $, respectively. The same gains are set for all the CAVs because we consider a homogeneous scenario. 

As mentioned earlier, we consider second-order vehicle dynamics as shown in~\eqref{equ:StateSpaceEquation}. We define the car-following errors as the new vehicle state.
\begin{equation}
	\bar{\boldsymbol{x}}_{i}=\left[
	\begin{array}{c}
		\bar{x}_{i, 1} \\
		\bar{x}_{i, 2} \\
	\end{array}\right]=\left[
	\begin{array}{c}
		p_{0}-p_{i}-D_\mathrm{des}\left(d_{0}-d_{i}\right) \\
		v_{i}-v_{\mathrm{p}} \\
	\end{array}\right], i \in \mathbb{N}^{+}.
\end{equation}
The vehicle input remains the same,~\ie, $ \bar{u}_{i} = u_{i} $.

Therefore, the car-following vehicle dynamic model is
\begin{equation}
	\dot{\bar{\boldsymbol{x}}}_{i}=\boldsymbol{A} \bar{\boldsymbol{x}}_{i}+\boldsymbol{B} \bar{u}_{i}, \, i \in \mathbb{N}^{+}.
\end{equation}

The linear feedback controller is simplified to
\begin{equation}
	\bar{u}_{i}=-k_{p} \sum_{j}\left(l_{i j}+q_{i j}\right) \bar{x}_{j, 1}-k_{v} \sum_{j}\left(l_{i j}+q_{i j}\right) \bar{x}_{j, 2},\, j \in \mathbb{I}_{i}.
\end{equation}

Defining $ \boldsymbol{k} = \left[k_{p},k_{i}\right]^{T} $, we have
\begin{equation}
	\bar{u}_{i}=-\sum_{j} \left(l_{ij}+q_{ij}\right) \boldsymbol{k}^{T}\bar{\boldsymbol{x}}_{i}, \, i \in \mathbb{N}^{+}.
\end{equation}

\begin{table}[!t]
	\renewcommand{\arraystretch}{1.3}
	\centering
	\begin{tabular}{c|c|c|c}
		\hline
		{\textbf{Types}} & {\textbf{Parameter}} &  {\textbf{Symbol}} & {\textbf{Value}} \\
		\hline
		\multirow{3}{*}{\textbf{Simulation}} & Simulation step & - &  $ 0.1 \mathrm{s} $ \\
		& Lane changing zone length & $ L_{\mathrm{LCZ}} $ & $ 500 \mathrm{m} $ \\
		& Car-following zone length & $ L_{\mathrm{CFZ}} $ & $ 500 \mathrm{m} $ \\
		\hline
		\multirow{3}{*}{\textbf{RCS}} & Group size & - &  $ 3 $ \\
		& One-step time & - & $ 4 \mathrm{s}$ \\
		& Lane changing time & - & $ 3 \mathrm{s}$ \\
		\hline
		\multirow{4}{*}{\textbf{Controller}} & Feedback gain of distance error & $ k_{p} $ & 0.1\\
		& Feedback gain of velocity error & $ k_{i} $ & $ 0.3 $\\
		& Steady platoon velocity & $ v_{\mathrm{p}} $ & $ 10 \mathrm{m/s} $\\
		& Steady car-following distance & $ d_\mathrm{f} $ & $ 30 \mathrm{m} $\\
		\hline
		\multirow{4}{*}{\textbf{Constraints}} & Maximum acceleration & $ a_{\max} $ & $5 \mathrm{m/s}^{2} $ \\
		& Minimum acceleration & $ a_{\min} $ &  $ -6 \mathrm{m/s}^{2} $ \\
		& Maximum velocity  & $ v_{\max} $ & $ 15 \mathrm{m/s} $ \\
		& Minimum velocity & $ v_{\min} $ & $ 0 \mathrm{m/s} $  \\
		\hline
	\end{tabular}
	\caption{Parameters}
	\label{tab:Parameters}
\end{table}

\section{Simulations}
\label{Sec:Simulation}
\subsection{Simulation Environment}
The traffic simulation was conducted in SUMO, which is widely used in traffic research\cite{lopez2018microscopic}. The simulation was run on an Intel Core i7-7700 @3.6 GHz processor. The intersection scenario and lane direction settings are the same as shown in Fig.~\ref{fig:Intersection}, and the overall length of the lane changing and car-following zones is $ 1000 \, m $. The CAV arrival is assumed to be a Poisson distributed flow, given by
\begin{equation}
	\label{equ:Poisson}
	P(X=k)=\frac{\lambda^{k}}{k!} e^{-\lambda}, k=0,1,\cdots,
\end{equation}
where $ X $ represents the arrival of the vehicle at the control zone. $ \lambda $ is the expected value as well as the variance of the Poisson distribution. Other key simulation parameters are listed in Table~\ref{tab:Parameters}.



\subsection{Evaluation Indexes}
As mentioned earlier, we aim to propose a cooperation method to improve both traffic safety and efficiency,~\ie, obtain the collision-free optimal CAV passing order. We assign $ t_{i}^{\mathrm{in}} $ as the time step when vehicle $ i $ enters the lane changing zone and $ t_{i}^{\mathrm{out}} $ as the time step when it arrives at the intersection. Several evaluation indexes have been proposed to measure the scheduling performance. In this study, we chose the evacuation time as the primary optimization target, as shown in Definition~\ref{def:EvacuationTime}. The evacuation time is related to the spanning tree depth, which is described in Section~\ref{Sec:StageTwo}.

\begin{definition}[Evacuation Time]
	\label{def:EvacuationTime}
	The evacuation time of $ N $ CAVs is defined as the time when the last CAV reaches the stop line; it is expressed as
	\begin{equation}
		\label{equ:EvacuationTime}
		t_\mathrm{evc} = \max t_{i}^{\mathrm{out}}, i \le N, i \in \mathbb{N}^+.
	\end{equation}
	Considering $ N $ incoming CAVs, $ t_\mathrm{evc} $ represents the arrival time of the last CAV at the stop line. For $ N $ CAVs, smaller evacuation times indicate that these CAVs pass through the intersection in a shorter time. Thus, it demonstrates the overall traffic efficiency performance,~\ie, the overall benefits to the CAVs. 
\end{definition}

In addition, the vehicle travels through the control zone in $ t_{i}^{\mathrm{out}} - t_{i}^{\mathrm{in}} $ time, whereas it travels through it under the free driving condition in $ {L_{\mathrm{ctrl}}}/{v_{\max}} $ time. Accordingly, the ATTD is selected to measure the average traffic efficiency of the vehicles, as described in Definition ~\ref{def:ATTD}. 
\begin{definition}[Average Travel Time Delay]
	\label{def:ATTD}
	The ATTD is designed to evaluate the average traffic efficiency of $ N $ CAVs. It is expressed as
	\begin{equation}
		\label{equ:TravelTimeDelay}
		t_\mathrm{ATTD} = \frac{1}{N} \sum_{i = 1}^{N} \left(t_{i}^{\mathrm{out}} - t_{i}^{\mathrm{in}}-\frac{L_{\mathrm{ctrl}}}{v_{\max}}\right),
	\end{equation}
	where $ L_{\mathrm{ctrl}} $ is the length of the control zone and $ t_\mathrm{ATTD} $ represents the average travel delay of the CAVs. The travel time of every CAV is considered in $ t_\mathrm{ATTD} $; therefore, it denotes the individual benefits for $ N $ CAVs, which is the secondary optimization target of traffic efficiency.   
\end{definition}

\begin{figure}[!t]
	\centering
	\includegraphics[width=\linewidth]{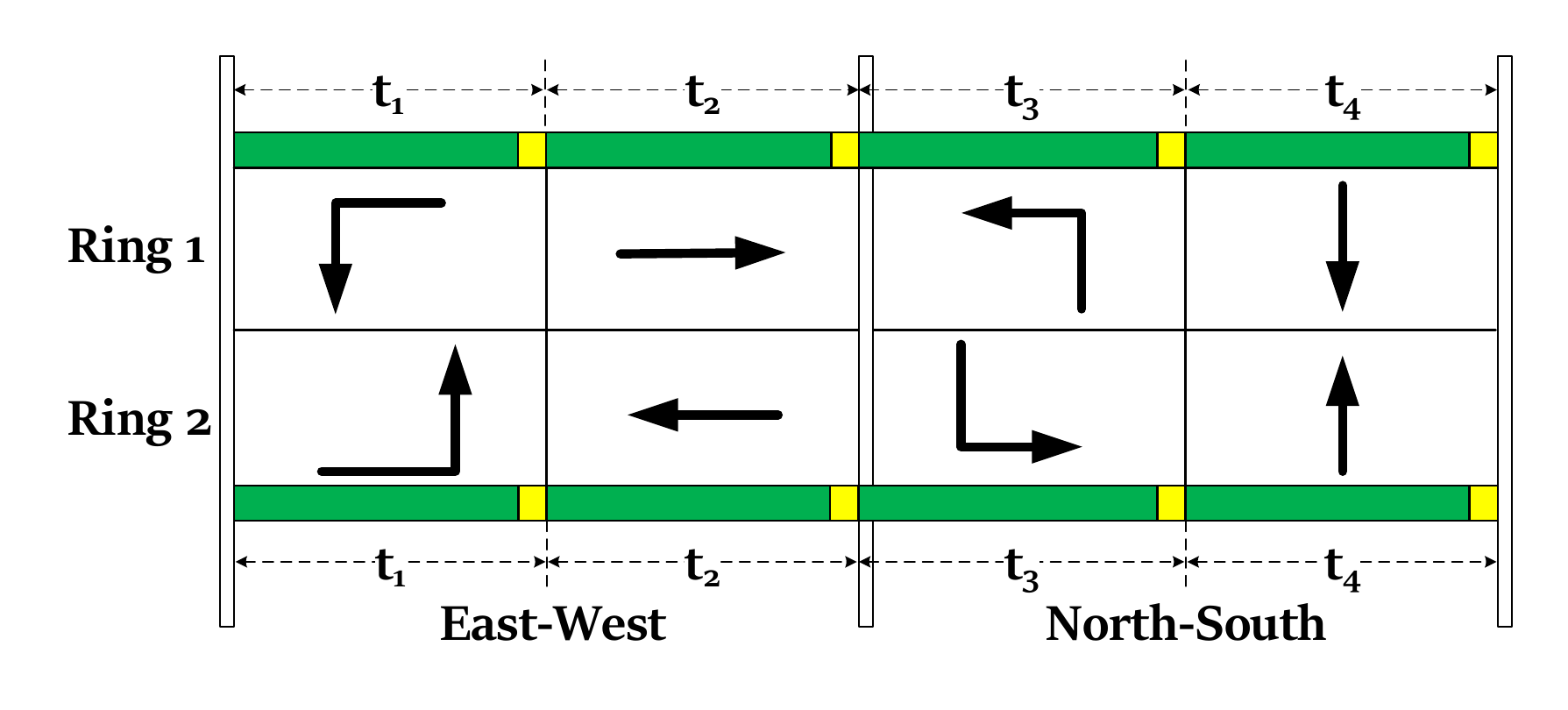}
	\caption{Dual ring control in constant traffic SPAT. The green line indicates that the corresponding traffic light is set to green, and the yellow square denotes the clearance time.}
	\label{fig:ConstantSPAT}
\end{figure}

\begin{figure*}[t]
	\centering
	\subcaptionbox{Evacuation time results comparison. \label{fig:VehicleNumber_Duration}}
	{\includegraphics[width=\linewidth]{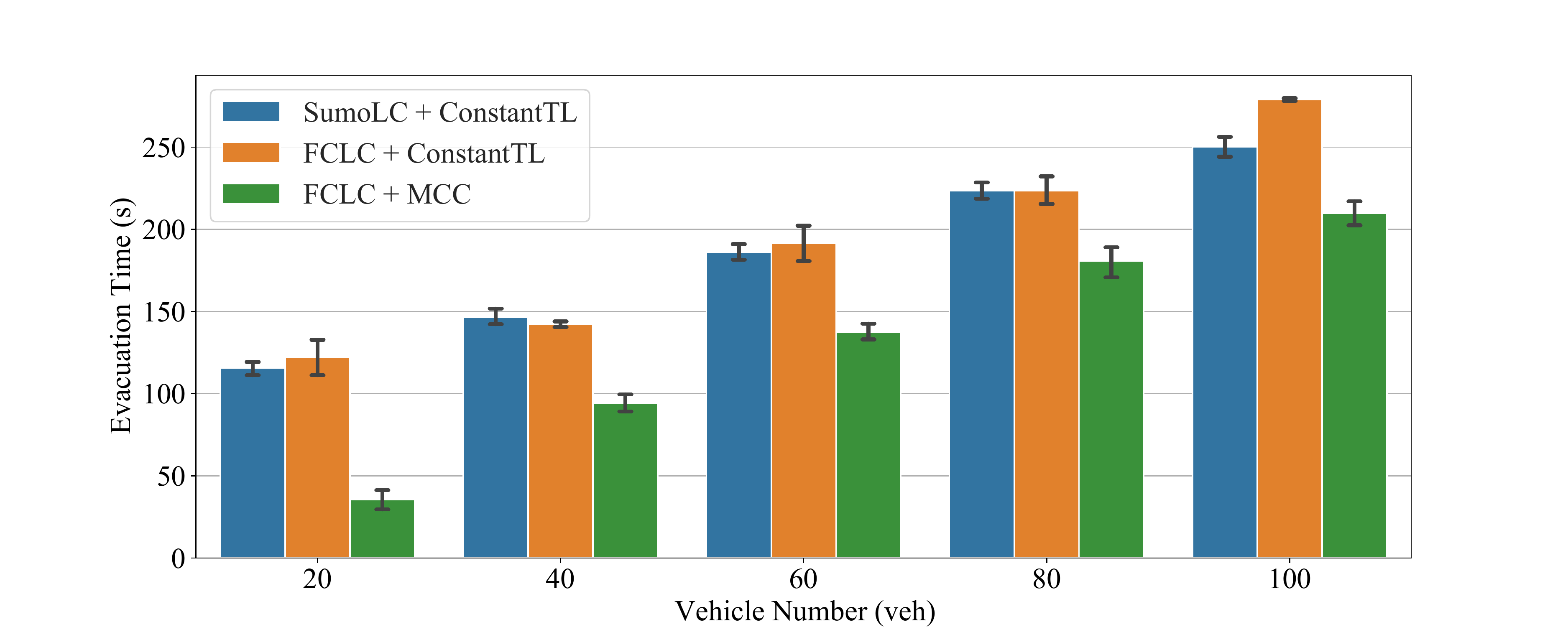}}
	\subcaptionbox{ATTD results comparison. \label{fig:VehicleNumber_ATTD}}
	{\includegraphics[width=\linewidth]{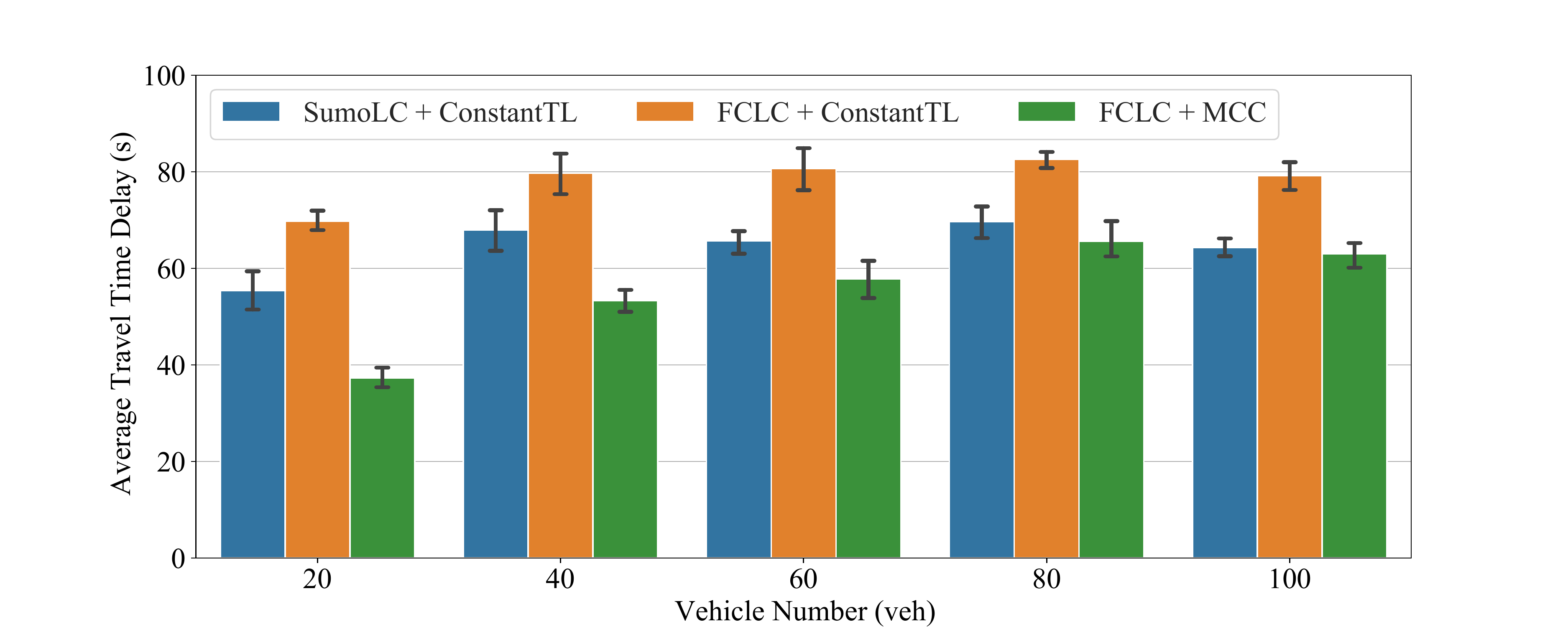}}
	\caption{Comparison of the traffic evacuation times and ATTDs of the algorithms for different numbers of input vehicles.}
	\label{fig:VehicleNumber}
\end{figure*}

As mentioned earlier, the control zone is divided into the lane changing and car-following zones. Although we primarily focus on the feasible path planning solution in the lane changing zone as described in Section~\ref{Sec:StageOne}, here, we continue to measure the algorithm performance over the complete length of the control zone.

\subsection{Benchmark algorithm}
In~\cite{chen2021conflict}, we have proved the optimality of the MCC algorithm in scheduling CAVs. This study further considers lane changing behavior; therefore, we compare our algorithm with the constant traffic SPAT method. The lanes leading to the intersection are divided into the lane changing and car-following zones, as shown in Fig.~\ref{fig:RoadSegmentation}; therefore, two algorithms are applied to each zone.


\subsubsection{Lane changing zone}
The formation control method is used as described in Section~\ref{Sec:StageOne} to obtain a parallel structure formation; therefore, the proposed lane changing zone algorithm is named as the formation-control lane changing (FCLC). As mentioned earlier, SUMO provides a lane changing model~\cite{erdmann2015sumo} to simulate the driver lane changing behavior. We use the default lane changing model as the benchmark algorithm, where the lane change timing is determined by the lane changing model rather than the formation control methods. Namely, the default lane changing model is called Sumo lane changing (SumoLC).

\begin{figure*}[t]
	\centering
	\subcaptionbox{Evacuation time results comparison. \label{fig:Lambda_Duration}}
	{\includegraphics[width=\linewidth]{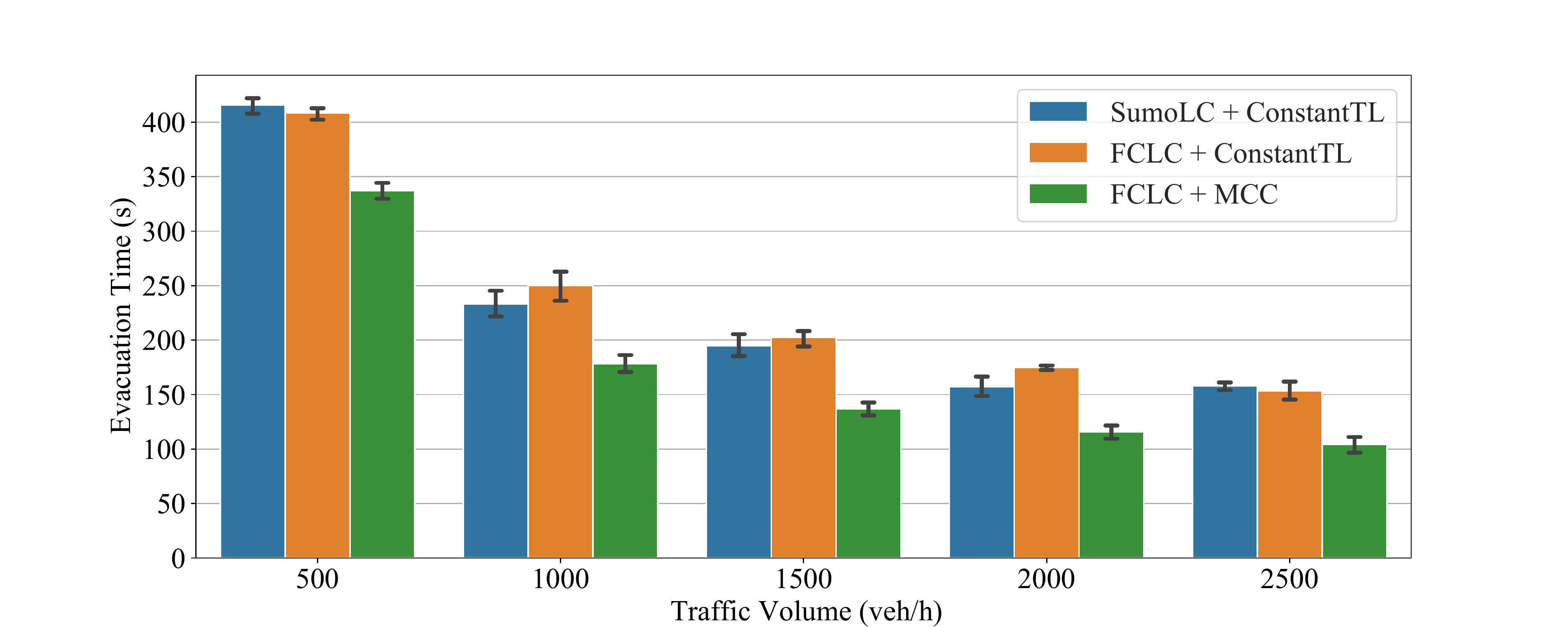}}
	\subcaptionbox{ATTD results comparison. \label{fig:Lambda_ATTD}}
	{\includegraphics[width=\linewidth]{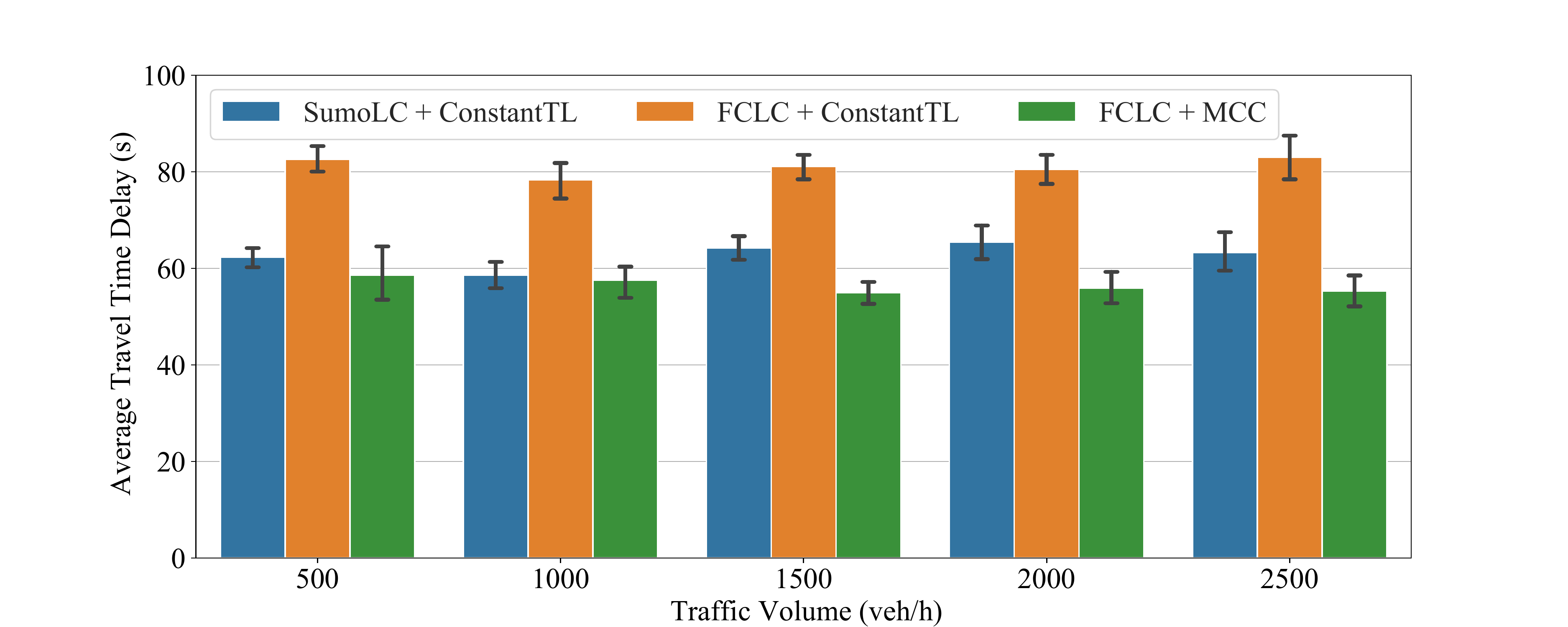}}
	\caption{Comparison of the traffic evacuation times and ATTDs of algorithms for different traffic volumes.}
	\label{fig:Lambda}
\end{figure*}

\subsubsection{Car-following zone}
As described in Section~\ref{Sec:StageTwo}, an MCC formulation is used to solve the scheduling problem; therefore, the MCC is referred to as the proposed algorithm. Constant traffic light (ConstantTL) is set as a dual ring control~\cite{xu2017v2i}, as shown in Fig~\ref{fig:ConstantSPAT}. In $ t_{1} $, the traffic light is set to green for vehicles turning left from either the east or the west. The yellow square represents the clearance time when the lights turn yellow. $ t_{2} $ represents the vehicles going straight from either the east or the west, and so on. A time of $ 35 \, \mathrm{s} $ is set for each phase, implying that
\begin{equation}
	t_{1} = t_{2} = t_{3} = t_{4} = 35 \, \mathrm{s},
\end{equation}
with $ 5 \, \mathrm{s} $ as the clearance time.

Theoretically, the combinations of two lane changing algorithms and two scheduling algorithms lead to four algorithms. However, SumoLC algorithm has its limitation that deadlock problem may occur if approaching lane is not long enough~\cite{erdmann2015sumo}. In our simulation, when vehicle number exceeds $ 40 $ or traffic volume reaches $ 1500 \, \mathrm{veh/h} $, lane changing failure occurs in one third of the results in SumoLC+MCC method. In comparison, deadlock does not happen in SumoLC+ConstantTL method because the control zone $ L_{\mathrm{ctrl}} $ is long enough. Moreover, FCLC guarantees the collision-free paths set is feasible. Similar results have been observed in~\cite{he2018erasing, cai2021formation}. If lane changing fails, CAVs cannot pass the intersection since there exist directional limitations on each lane. Therefore, we excluded this method. The comparison of the other three algorithms still provides convincing proof on the effectiveness of the proposed algorithm.

\subsection{Simulation Results for Different Numbers of Input Vehicles}

The first simulation was conducted for different numbers of vehicles. The scheduling problem is highly related to the distribution of incoming vehicles; therefore, ten iterations of the simulation were conducted for each number of vehicles and for each algorithm. Notably, identical vehicle distributions were applied for the algorithms,~\ie, ten sets of vehicle distributions were randomly generated for each number of vehicles and used as the same input for the simulated algorithms. Typically, the number of CAVs in one intersection is less than $ 100 $; thus, we set the number of vehicles to range from $ 20 $ to $ 100 $. The traffic volume was set as $ 2000 \, \mathrm{veh/lane} $.

The simulation result is shown in Fig.~\ref{fig:VehicleNumber}. The bar plot shows the average value of the ten sets of simulations. The standard deviations of each set of results are also provided. Generally, both the evacuation time and the ATTD increase with the number of vehicles because of vehicle queuing. It is apparent that the MCC in the lane changing zone significantly improves the traffic efficiency. The evacuation time is the primary optimization target in the MCC scheduling; therefore, the MCC algorithm reduces the evacuation time by $ 15.8\% $ to $ 67.9\% $ for $ 20 $ to $ 100 $ vehicles. 

Regarding the ATTD perspective, it is evident that the FCLC negatively influences the ATTD. This is because when SumoLC is applied, CAVs are set to the maximum velocity $ v_{\max} = 15 \, \mathrm{m/s} $ if there are no other vehicles around. However, the FCLC aims to develop a steady parallel formation with a platoon velocity $ v_{\mathrm{p}} = 10 \, \mathrm{m/s} $. Therefore, the total evacuation time is unaffected when vehicle number is less than $ 80 $ because the intersection is not oversaturated; however, the ATTD is higher because of the delay in the lane changing zone. On the other side, the formation control in FCLC helps to generate a collision-free and feasible path planning set to MCC. Moreover, the steady parallel formation is convenient for MCC to adjust the car following distance to desired virtual platoon positions. MCC algorithm reduces $ 2.2 \% $ to $ 32.6 \% $ of the ATTD. Moreover, as the vehicle number increases, it has less improvement on ATTD since ATTD is not the priority concern in optimization.




\subsection{Simulation Results for Different Traffic Volumes}

The second simulation was conducted for different traffic volumes. Similar to the previous simulation, ten randomly generated repetitions of the simulation were conducted for each traffic volume. The number of vehicles was set to $ 50 $, and the traffic volume varied from $ 500 \, \mathrm{veh/h} $ to $ 2500 \, \mathrm{veh/h} $.


Generally, as the number of vehicles is constant, the evacuation time is decreased as the traffic volume increases. The MCC exhibits noticeable improvement in the traffic efficiency, which reduces the evacuation time by $ 18.8\% $ to $ 33.2\% $ and ATTD by $ 1.9\% $ to $ 14.5\% $. The comparison between the SumoLC and FCLC results agrees with the analysis of the simulation results for different numbers of vehicles. FCLC delays the travel time in lane changing zone because of velocity settings. FCLC+MCC provides a feasible solution with noticeable improvement in the traffic efficiency.

\section{Conclusions}
\label{Sec:Conclusion}

In this paper, a two-stage cooperation framework is proposed to improve vehicle safety and traffic efficiency at intersections where lane changing is permitted. In the first stage, we design an iterative method to solve the multi-vehicle target assignment and path planning problem. In contrast to existing single-vehicle lane changing algorithms, the deadlock problem is solved in our method and a feasible collision-free path planning is ensured. In the second stage, a graph-based method is proposed to schedule the CAV arrival time. A heuristic algorithm is established to solve the problem with a low computational burden. Traffic simulations verified the effectiveness of the proposed algorithm.

Future research directions include studying the relationship of the two decoupled stages further. In this study, the optimal arrival plan is calculated when the CAV arrives at the car-following zone. However, it can also be calculated when the CAV enters the lane changing zone. Thus, in the lane changing zone, a CAV can drive to its target position in virtual platoon. Another interesting topic involves considering HDVs in this scenario because lane changing would be permitted. There has been little research on intersections with mixed traffic where lane changing is permitted.

\ifCLASSOPTIONcaptionsoff
  \newpage
\fi



\bibliographystyle{IEEEtran}
\bibliography{IEEEabrv,mybibfile}
\end{document}